\input luatex-or-tex.mac

\citenumber

\null

\vskip 0.8 true cm

\centerline{\bf A note on deformation argument for $L^2$ constraint problems}

\bigskip

\centerline{Norihisa Ikoma\footnote{${}^*$}{The first author is partially supported
by JSPS KAKENHI Grant No. JP16K17623, JP17H02851.}
and 
Kazunaga Tanaka\footnote{${}^{**}$}{The second author is partially supported
by JSPS KAKENHI Grant No. JP17H02855, JP16K13771, JP26247014, JP18KK0073.\m
{\tt MSC2010: 35J20, 35J50, 35Q55, 58E05}
}
}

\medskip

\settabs 27 \columns
\+&${}^*$& Department of Mathematics, Faculty of Science and Technology, Keio University\cr
\+&& Yokohama, Kanagawa 223-8522, Japan\cr
\+&& email: {\tt \raw[ikoma@math.keio.ac.jp]}\cr
\+&${}^{**}$& Department of Mathematics, School of Science and Engineering, Waseda University\cr
\+&& 3-4-1 Ohkubo, Shinjuku-ku, Tokyo 169-8555, Japan\cr
\+&& email: {\tt \raw[kazunaga@waseda.jp]}\cr

\bigskip


\input macro.mac


{\narrower
\noindent {\bf Abstract.}
We study the existence of $L^2$ normalized solutions for nonlinear Schr\"odinger equations and 
systems.   Under new Palais-Smale type conditions we develop new deformation arguments 
for the constraint functional on $S_m=\{ u; \, \int_{\R^N}\abs u^2=m\}$ or $\SS$.
As applications, we give other proofs to the results of [\cite[J:20], \cite[BdV:6], \cite[BS1:7]].
As to  the results of [\cite[J:20], \cite[BdV:6]], our deformation result enables us to apply
the genus theory directly to the corresponding functional to obtain infinitely many solutions.
As to the result [\cite[BS1:7]], via our deformation result we can show the existence of 
vector solution without using constraint related to the Pohozaev identity.

}

\medskip

\BS{\label[Section:1]. Introduction}
In this paper we develop a new deformation argument for $L^2$-constraint problems and as applications
we study the existence and multiplicity of standing waves $\Phi(t,x)=e^{i\lambda t}u(x)$
of nonlinear Schr\"odinger equations
    $$  -i\partial_t\Phi =\Delta \Phi +g(\Phi) \quad (t,x)\in \R\times\R^N
    $$
and those $\Phi_i(t,x)=e^{i\lambda_i t}u_i(x)$ ($i=1,2$) of nonlinear Schr\"odinger systems
    $$  \left\{\eqalign{
        &-i\partial_t\Phi_1 = \Delta\Phi_1 +\mu_1\abs{\Phi_1}^2\Phi_1 +\beta\abs{\Phi_2}^2\Phi_1,
                \quad  (t,x)\in \R\times\R^3, \cr
        &-i\partial_t\Phi_2 = \Delta\Phi_2 +\mu_2\abs{\Phi_2}^2\Phi_2 +\beta\abs{\Phi_1}^2\Phi_2,
                \quad  (t,x)\in \R\times\R^3.\cr}
        \right.
    $$
It is easily seen that $u$ ($(u_1,u_2)$ respectively) are solutions of 
    $$  \eqalign{
        &\quad -\Delta u +\lambda u= g(u) \quad \hbox{in}\ \R^N,   \cr
        &\left(
        \left\{\eqalign{
            -&\Delta u_1+\lambda_1u_1 = \mu_1u_1^3 + \beta u_1u_2^2 \quad \hbox{in}\ \R^3, \cr
            -&\Delta u_2+\lambda_2u_2 = \mu_2u_2^3 + \beta u_1^2u_2 \quad \hbox{in}\ \R^3, \cr}
        \right.
        \quad \hbox{respectively} \right). \cr}
    $$
We study the existence of {\it normalized solutions}, that is, for given $m$ and $(m_1,m_2)$, we try
to find solutions of \ref[1.1] and \ref[1.2] below.

When we take variational approaches for nonlinear elliptic problems, deformation arguments play 
essential roles and they enable us to apply minimax methods 
and other topological tools (e.g. genus theories, symmetric mountain pass theorems etc.), which
give us existence and multiplicity results.
However deformation arguments are not developed well for $L^2$-constraint problems and
slightly different approaches, which are in the spirit of Ekeland' variational principle, are taken.

More precisely, we consider nonlinear Schr\"odinger equations:
    $$  \left\{\eqalign{
        -&\Delta u +\lambda u= g(u) \quad \hbox{in}\ \R^N, \cr
        &\intRN\abs u^2=m, \cr}\right.
        \eqno\label[1.1]
    $$
where $m>0$ is a given constant and $g(\xi)\in C(\R,\R)$ is a function with $L^2$ super-critical 
growth and Sobolev sub-critical growth (See (g1)--(g2) below).  In [\cite[J:20]], Jeanjean
finds a solution of \ref[1.1] as a critical point of the following functional
    $$  I(u)=\half\intRN\abs{\nabla u}^2 - \intRN G(u):\, S_m\to \R,
    $$
where $G(\xi)=\int_0^\xi g(\tau)\, d\tau$ and $S_m=\{ u\in\E;\, \intRN\abs u^2=m\}$.
He also introduces an augmented functional
    $$  J(\theta,u)=I(e^{{N\over 2}\theta} u(e^\theta u)) 
        = \half e^{2\theta}\intRN\abs{\nabla u}^2 - e^{-N\theta} \intRN G(e^{{N\over 2}\theta} u)
        :\, \R\times S_m\to \R.
    $$
A key step of his argument is to generate a sequence $\{ (\theta_j,u_j)\}_{j=1}^\infty$ such that
    $$  \theta_j\to 0, \ \partial_\theta J(\theta_j,u_j)\to 0, \ 
        \norm{\partial_u J(\theta_j,u_j)}_{(\E)^*} \to 0, \ J(\theta_j,u_j)\to b.
    $$
We remark that $\partial_\theta J(0,u)=0$ is related to the Pohozaev identity and the sequence
$\{ (\theta_j,u_j)\}_{j=1}^\infty$ is a Palais-Smale sequence for $J$ with an extra property 
$\partial_\theta J(\theta_j,u_j)\to 0$, which makes it possible to extract a convergent subsequence.
We also refer to [\cite[BdV:6]], in which Bartsch and de Valeriola show the existence of infinitely
many solutions of \ref[1.1] via \lq\lq fountain theorem'' under the assumption of oddness of $g(\xi)$.

\medskip

\proclaim Theorem \label[Theorem:1.1] {\rm ([\cite[J:20], \cite[BdV:6]])}.  
Assume $N\geq 2$, $m>0$ and
\itemitem{(g1)} $g(\xi):\, \R\to\R$ is continuous.
\itemitem{(g2)} There exist constants $\alpha$, $\beta\in\R$ with $2+{4\over N}<\alpha<\beta<2^*$ such
that
    $$  0 < \alpha G(\xi) \leq g(\xi)\xi \leq \beta G(\xi) \quad \hbox{for}\ \xi\in\R\setminus\{ 0\},
    $$
where $2^*={2N\over N-2}$ for $N\geq 3$, $2^*=\infty$ for $N=2$.

{\sl\noindent
Then we have
\item{(i)} \ref[1.1] has at least one solution.
\item{(ii)}  In addition to (g1)--(g2), assume $g(-\xi)=-g(\xi)$ for $\xi\in\R$.
Then \ref[1.1] has infinitely many solutions.

}

\medskip

\noindent
We also note that a similar approach was taken for nonlinear scalar field equations in Hirata, Ikoma
and Tanaka [\cite[HIT:17]] successfully and they gave another proof of the results in [\cite[BeL1:9], 
\cite[BeL2:10], \cite[BeGK:11]].  We also refer to 
[\cite[AdAP:3], \cite[BT:12], \cite[CT:13], \cite[FISJ:14], \cite[I1:18], \cite[I2:19], \cite[MVS:24]] for 
other applications of generation of Palais-Smale sequences with an extra property.

Recently Bartsch and Soave [\cite[BS1:7]] (c.f. [\cite[BS2:8]]) study an $L^2$-constraint problem using {\it natural constraint}.
They consider the following systems of nonlinear Schr\"odinger equations:
    $$  \left\{\eqalign{
        -&\Delta u_1+\lambda_1u_1 = \mu_1u_1^3 + \beta u_1u_2^2 \quad \hbox{in}\ \R^3, \cr
        -&\Delta u_2+\lambda_2u_2 = \mu_2u_2^3 + \beta u_1^2u_2 \quad \hbox{in}\ \R^3, \cr
        &\int_{\R^3} \abs{u_1}^2 =m_1, \quad \int_{\R^3} \abs{u_2}^2 =m_2, \cr}\right.
                                                                    \eqno\label[1.2]
    $$
where $\mu_1$, $\mu_2>0$, $\beta<0$, $m_1$, $m_2>0$ are given constants and $\lambda_1$, $\lambda_2\in\R$
are unknown Lagrange multipliers.  We remark that they deal with the focusing-repulsive case 
$\mu_1$, $\mu_2>0$, $\beta<0$ and solutions of \ref[1.2] are characterized as critical points of
    $$  I_*(u_1,u_2) = \half\int_{\R^3} \abs{\nabla u_1}^2+\abs{\nabla u_2}^2
        -{1\over 4}\int_{\R^3} \mu_1 u_{1}^4+\mu_2 u_{2}^4 +2\beta u_1^2u_2^2
        : \SS\to\R.
    $$
To find critical points of $I_*(u_1,u_2)$, they introduced a {\it natural constraint}:
    $$  \calP=\{ (u_1,u_2)\in \SS;\, P_*(u_1,u_2)=0\},
    $$
and
    $$  P_*(u_1,u_2)= \int_{\R^3}\abs{\nabla u_1}^2 + \abs{\nabla u_2}^2
            - {3\over 4}\int_{\R^3}\mu_1u_{1}^4 +\mu_2u_{2}^4 +2\beta u_1^2u_2^2.
    $$
Here $P_*$ is a functional related to the Pohozaev identity and, using the result of Ghoussoub [\cite[G:15]],
they showed

{\parindent=1.5\parindent

\item{(P-i)} $\calP$ is a $C^1$ submanifold of $\SS$;

\item{(P-ii)} For a suitable minimax value $c$, there exists a Palais-Smale sequence \m
$\{(v_{1n}, v_{2n})\}_{n=1}^\infty\subset\calP$ for $I_*$ at level $c$.

}

\noindent
From these properties they show the existence of critical value of $I_*$.  Thus they obtain

\medskip

\proclaim Theorem \label[Theorem:1.2]  ([\cite[BS1:7]]).
Let $N=3$ and let $\mu_1$, $\mu_2>0$, $\beta<0$, $m_1$, $m_2>0$.  Then \ref[1.1] has a solution
$(\lambda_1,\lambda_2,u_1,u_2)$ with $\lambda_1$, $\lambda_2>0$ and $u_1$, $u_2$ are positive in $\R^3$
and radially symmetric.

\medskip

\noindent
We remark that in [\cite[BS1:7]] they also develop natural constraint approach for \ref[1.1] under
additional condition.  See Remark \ref[Remark:2.4].
We also refer to [\cite[AKY:1], \cite[AN:2], \cite[AP:4], \cite[LM:22], \cite[Sh:27]] for application of 
the Pohozaev manifolds for nonlinear Schr\"odinger equations (without $L^2$-constraint).  

\medskip

\noindent
In this paper, we introduce new deformation approaches for $I(u)\in C^1(S_m,\R)$ and \m
$I_*(u_1,u_2)\in C^1(\SS,\R)$.  The main
difficulty to deal with the corresponding functionals to \ref[1.1] on $S_m$ or to \ref[1.2]
on $\SS$ is the lack of the Palais-Smale condition, that is, it is difficult to verify
the Palais-Smale condition and thus the usual deformation argument does not work directly 
for these problems.  
To overcome this difficulty we introduce a new type of Palais-Smale condition $(PSP)_c$ for 
$c\in \R$.  We just state it  for the functional $I_*(u_1,u_2):\, \SS\to\R$ corresponding to 
\ref[1.2]. 

\smallskip

\itemitem{$(PSP)_c$} If a sequence $\{(u_{1n}, u_{2n})\}_{n=1}^\infty\subset\SS$ satisfies
    $$  \eqalign{
        &I_*(u_{1n},u_{2n})\to c, \cr
        &\norm{dI_*(u_{1n},u_{2n})}_{T_{(u_{1n},u_{2n})}^*(\SS)} \to 0,\cr
        &P_*(u_{1n},u_{2n})\to 0, \cr}
    $$
then $\{(u_{1n}, u_{2n})\}_{n=1}^\infty\subset\SS$ has a strongly convergent subsequence.  
In particular, $c$ is a critical point value of $I_*$.

\smallskip

\noindent
Our deformation result for $I_*\in C^1(\SS,\R)$ is

\proclaim Proposition \label[Proposition:1.3].
For $c\in \R$, suppose that $(PSP)_c$ holds.  Then for any $\overline\epsilon>0$ and any neighborhood
$O$ of $K_c=\{ (u_1,u_2)\in \SS;\, I_*(u_1,u_2)=c,\, dI_*(u_1,u_2)=0, \}$, there exist
$\epsilon\in (0,\overline\epsilon)$ and $\eta\in C([0,1]\times\SS,\SS)$ such that
\item{(i)} $\eta(0,u_1,u_2)=(u_1,u_2)$ for $(u_1,u_2)\in\SS$.
\item{(ii)} $\eta(t,u_1,u_2)=(u_1,u_2)$ for $t\in [0,1]$ if $I_*(u_1,u_2)\leq c-\overline\epsilon$.
\item{(iii)} $t\mapsto I_*(\eta(t,u_1,u_2))$ is non-increasing for all $(u_1,u_2)\in \SS$.
\item{(iv)} $\eta(1,[I_*\leq c+\epsilon]\setminus O) \subset [I_*\leq c-\epsilon]$, 
$\eta(1,[I_*\leq c+\epsilon]) \subset [I_*\leq c-\epsilon]\cup O$, \m
where
    $$  [I_*\leq c] =\{(u_1,u_2)\in \SS;\, I_*(u_1,u_2)\leq c\} \quad \hbox{for}\ c\in\R.
    $$

\noindent
Proposition \ref[Proposition:1.3] enables us to apply minimax methods to find critical points of $I_*$.
We note that the $(PSP)_c$ condition is weaker than the Palais-Smale condition.
In the approaches in [\cite[BS1:7], \cite[BS2:8]], they need to introduce a natural constraint.
In our approach we emphasis that 
we don't need to introduce any constraint to $I_*$.
Our Proposition \ref[Proposition:1.3] will be derived from Proposition \ref[Proposition:4.5], in which the deformation result is obtained
in a general setting.  See Remark \ref[Remark:4.6] for other merits of our approach.

To show the deformation result in Proposition \ref[Proposition:1.3] (also in Proposition \ref[Proposition:4.5]),
we use an idea from Hirata and Tanaka [\cite[HT:16]].
Here we give a heuristic explanation in the setting of Proposition \ref[Proposition:1.3].  
In our deformation argument, the following $\R$-action $\Phi_\theta$ on $\SS$ is important:
    $$  \Phi_\theta(u_1,u_2)(x)=(e^{{3\over 2}\theta}u_1(e^\theta x),e^{{3\over 2}\theta}u_2(e^\theta x))
        \quad \hbox{for}\ \theta\in\R \ \hbox{and}\ (u_1,u_2)\in \SS.
    $$
We note that $\SS$ is invariant under $\Phi_\theta$ and
    $$  P_*(u_1,u_2)={d\over d\theta}\Big|_{\theta=0} I_*(\Phi_\theta(u_1,u_2)).
    $$
The equation 
$P_*(u_1,u_2)=0$ is related to the Pohozaev identity and it is natural to expect that the Pohozaev
identity holds for any solution and it is verified in Proposition \ref[Proposition:3.1].

Conversely if the Pohozaev identity does not hold, i.e., 
$P_*(u_1,u_2)>0$ ($P_*(u_1,u_2)<0$ respectively), for $(v_1,v_2)$ in a small neighborhood
of $(u_1,u_2)$ and for small $\delta>0$, the following map gives a continuous deformation on
$\SS$
    $$  [0,\delta]\to \SS;\ t\mapsto \Phi_{-t}(v_1,v_2) \quad
            (\Phi_t(v_1,v_2)\ \hbox{respectively}),     \eqno\label[1.3]
    $$
along which $I_*$ decreases.  Thus, it seems that critical points $(u_1,u_2)\in\SS$ of $I_*$
with $P_*(u_1,u_2)\not=0$ don't affect the topology of the level set $[I_*\leq c]$ much.
We note that the flow \ref[1.3] is continuous
but not of class $C^1$ and it seems that such a deformation cannot be obtained by the standard
deformation flow.
To justify such an observation, we argue in augmented space $\R\times \SS$; we set
    $$  \eqalign{
        J_*(\theta,u_1,u_2) &= I_*(\Phi_\theta(u_1,u_2)) \cr
        &= \half e^{2\theta} \int_{\R^3} \abs{\nabla u_1}^2+\abs{\nabla u_2}^2
            -{1\over 4} e^{3\theta} \int_{\R^3} \mu_1 u_{1+}^4+\mu_2 u_{2+}^4 +2\beta u_1^2u_2^2. \cr}
    $$
First we construct a deformation flow $\weta(t,\theta,u_1,u_2):\, [0,1]\times \R\times\SS\to 
\R\times \SS$ for $J_*(\theta,u_1,u_2)$ on $\R\times\SS$ and second we define a desired flow
$\eta(t,\theta,u_1,u_2):\, [0,1]\times \SS\to \SS$ by
    $$  \eta(t,u_1,u_2)=\pi(\weta(t,0,u_1,u_2)),
    $$
where $\pi:\, \R\times\SS\to \SS$ is given by
    $$  \pi(\theta,u_1,u_2)=\Phi_\theta(u_1,u_2).
    $$
We note that a map
    $$  [0,\delta]\to \R\times \SS;\, t\mapsto (-t,v_1,v_2) \qquad ((t,v_1,v_2) \ \hbox{respectively})
    $$
corresponds to \ref[1.3] and our deformation flow on $\SS$ is realized as
a composition $\pi$ and $C^1$-deformation $\weta$ in $\R\times\SS$.

We also note that in [\cite[HT:16]] we develop related deformation arguments for nonlinear
scalar field equations:
    $$  -\Delta u =g(u) \quad \hbox{in}\ \R^N, \quad u\in\E     \eqno\label[1.4]
    $$
and for $L^2$-constraint problems:
    $$  \left\{\eqalign{
        -&\Delta u +\lambda u=g(u) \quad \hbox{in}\ \R^N, \cr
        &\intRN \abs u^2 =m.  \cr } \right.         \eqno\label[1.5]
    $$
We study \ref[1.5] for $L^2$-subcritical nonlinearities in [\cite[HT:16]]. 
We also refer to [\cite[JL:21], \cite[Si:28]] for the studies of \ref[1.5].

Solutions of \ref[1.4] (\ref[1.5] respectively) can be characterized as critical points of
    $$  \eqalignno{
        &u\mapsto \half\intRN \abs{\nabla u}^2 -\intRN G(u);\, \E\to\R,     &\label[1.6]\cr
        &\Bigl( (\lambda,u)\mapsto \half\intRN \abs{\nabla u}^2 -\intRN G(u)
            +{\lambda\over 2}\bigl(\bigl(\intRN\abs u^2\bigr)-m\bigr); \cr
        &\qquad \qquad\qquad\qquad  \R\times\E\to\R,    
                    \quad \hbox{respectively}   \Bigr).                         &\label[1.7]\cr}
    $$
(We use the Lagrange formulation for \ref[1.5]).  

\noindent
We consider the corresponding augmented functional on $\R\times\E$ ($\R\times\R\times \E$ respectively)
and we succeeded to get deformation flows for \ref[1.6] and \ref[1.7], which enables us to 
give another proof to the results on [\cite[BeL1:9], \cite[BeGK:11]] on \ref[1.4] and to get a multiplicity
result for \ref[1.5].   
We also note that the deformation arguments are developed in the full spaces $\E$ and $\R\times\E$ in
[\cite[HT:16]].
Such a deformation theory can be developed also on the embedded manifolds (e.g. on $\SS$ for \ref[1.2])
and we give an abstract result in Section \ref[Section:4], which can be applied to \ref[1.1], \ref[1.2],
\ref[1.4] and \ref[1.5].

As applications of our new deformation argument, in this paper we deal with \m $L^2$-constraint problems 
for nonlinear Schr\"odinger equations \ref[1.1] and systems \ref[1.2].
In Section \ref[Section:2], we study \ref[1.1] and we give another proof of Theorem \ref[Theorem:1.1].  In particular, our new
deformation argument (Proposition \ref[Proposition:2.3]) enables us to apply the genus theory and symmetric mountain
pass theorem directly to the corresponding functional to obtain multiplicity result.

In Section \ref[Section:3], we study systems of nonlinear Schr\"odinger equations \ref[1.2], which is due to 
[\cite[BS1:7]],  and give a simpler proof
to Theorem \ref[Theorem:1.2].  We develop a new minimax methods for the functional $I_*(u_1,u_2):\, \SS\to\R$
and we give a minimax characterization of the solutions, which we believe of interest.  We also note
that scaling properties of $I_*$ and a Louiville type result, which is an extension of [\cite[BS1:7]],
play important roles in our minimax method.

Finally in Section \ref[Section:4], we give our deformation theory in an abstract setting.  It is used in 
Sections \ref[Section:2]--\ref[Section:3] and it also covers the results in [\cite[HT:16]].

\BS{\label[Section:2].  Single equations}
In this section we study the $L^2$-constraint problem for single equations and give other
approaches to the results of Jeanjean [\cite[J:20]] and Bartsch-de Valeriola [\cite[BdV:6]].
We also remark that results in Sections \ref[Subsection:2.1]--\ref[Subsection:2.2] are also 
important to study Schr\"odinger systems.

In what follows, we use notation:
    $$  \norm u_p = \left(\intRN \abs u^p\right)^{1/p} 
        \quad \hbox{for}\ p\in [1,\infty) \ \hbox{and}\ u\in L^p(\R^N).
    $$

\BSS{\label[Subsection:2.1]. Preliminaries}
We study the $L^2$-constraint problems for nonlinear Schr\"odinger equations:
    $$  \left\{\eqalign{
        -&\Delta u =g(u) -\lambda u \quad  \hbox{in}\ \R^N, \cr
        &\norm u_2^2 = m, \cr} \right.                      \eqno\label[2.1]
    $$
where $N\geq 2$, $m>0$ is a given constant, $g(\xi)$ is a given function satisfying the conditions (g1)--(g2)
in Theorem \ref[Theorem:1.1] and $\lambda\in\R$ is a unknown Lagrange multiplier.

Setting
    $$  \eqalign{
        &I(u) =\half\norm{\nabla u}_2^2 -\intRN G(u):\, \E \to\R, \cr
        &S_m = \{ u\in \E;\, \norm u_2^2=m\}, \cr}
    $$
solutions of \ref[2.1] can be characterized as critical points of $I\in C^1(S_m,\R)$.

For $u\in S_m$ and $t>0$, we set
    $$  u_t(x) = t^{N\over 2}u(tx).             \eqno\label[2.2]
    $$
We note that for $u\in\E$
    $$  \eqalignno{
        &\norm{u_t}_2^2 = \norm u_2^2,          &\label[2.3]\cr
        &\norm{\nabla u_t}_2^2 = t^2 \norm{\nabla u}_2^2,\cr
        &\norm{u_t}_r^r = t^{{r-2\over 2}N}\norm u_r^r.\cr}
    $$
In particular, we have $u_t\in S_m$ for all $u\in S_m$ and $t\in (0,\infty)$.

We also set
    $$  \eqalign{
        &P(u) = \norm{\nabla u}_2^2 -N\intRN \half g(u)u -G(u)\in C^1(\E,\R), \cr
        &\calP_m =\{ u\in S_m;\, P(u)=0\}.\cr}
    $$
We also note that any solution $u$ of \ref[2.1] satisfies the Pohozaev identity $P(u)=0$.
First we have

\proclaim Lemma \label[Lemma:2.1].
Assume (g1)--(g2).  Then we have
\item{(i)} For $u\in S_m$,
    $$  \eqalign{
        &I(u_t), \ P(u_t) \to -\infty \quad \hbox{as}\ t\to\infty,\cr
        &I(u_t), \ P(u_t) \to +0 \quad \ \hbox{as}\ t \to +0.\cr}
    $$
\item{(ii)} 
    $$  b_0 \equiv \inf_{u\in \calP_m} I(u)>0.          \eqno\label[2.4]
    $$

\noindent
To prove Lemma \ref[Lemma:2.1], we note that for some constants $C_1$, $C_2>0$,
    $$  \eqalignno{
        &C_1\min \{\abs \xi^\alpha, \abs \xi^\beta\} \leq G(\xi) \leq C_2\max \{\abs \xi^\alpha, \abs \xi^\beta\}
        \leq C_2(\abs \xi^\alpha + \abs \xi^\beta),         &\label[2.5]\cr
        &({\alpha\over 2}-1)G(\xi) \leq \half g(\xi) \xi -G(\xi) \leq ({\beta\over 2}-1)G(\xi)
        \leq ({\beta\over 2}-1)C_2(\abs \xi^\alpha + \abs \xi^\beta),
                                                &\label[2.6] \cr}
    $$
for all $\xi\in\R$.
These inequalities follow from (g2) easily.  We also use the following inequalities frequently, 
which follow from the Gagliard-Nirenberg inequality:
for all $u\in H^1(\R^N)$
    $$  \eqalignno{
        &\norm u_\alpha^\alpha \leq C_3 \norm{\nabla u}_2^{\alpha'}\norm u_2^{\alpha-\alpha'}, &\label[2.7]\cr
        &\norm u_\beta^\beta \leq C_3 \norm{\nabla u}_2^{\beta'}\norm u_2^{\beta-\beta'}, &\label[2.8]\cr}
    $$
where $\alpha'={\alpha-2\over 2}N\in(2,\alpha)$, $\beta'={\beta-2\over 2}N\in (2,\beta)$ and
$C_3>0$ is independent of $u$.

\medskip

\claim Proof of Lemma \ref[Lemma:2.1].
(i) For $u\in S_m$, we have from \ref[2.3] and \ref[2.5]
    $$  \eqalign{
        &\intRN G(u_t)
        \geq C_1\intRN \min\{ \abs{u_t}^\alpha, \abs{u_t}^\beta\} 
        = C_1 t^{{\alpha-2\over 2}N} \intRN \min\{ \abs u^\alpha, t^{{N\over 2}(\beta-\alpha)} \abs u^\beta\}\cr
        &\qquad \geq C_1 t^{{\alpha-2\over 2}N} \intRN \min\{ \abs u^\alpha, \abs u^\beta\}
            \quad \hbox{for}\ t\in [1,\infty), \cr
        &\intRN G(u_t)
        \leq C_2(\norm{u_t}_\alpha^\alpha + \norm{u_t}_\beta^\beta) 
        \leq C_2(t^{{\alpha-2\over 2}N}\norm{u}_\alpha^\alpha
            + t^{{\beta-2\over 2}N}\norm{u}_\beta^\beta)
            \quad \hbox{for}\ t\in (0,1].\cr}
    $$
Noting ${\beta-2\over 2}N>{\alpha-2\over 2}N >2$, we have (i). \m
(ii) For $u\in\calP_m$, by \ref[2.6]
    $$  \eqalign{
        0 &= \norm{\nabla u}_2^2 - N\intRN \half g(u)u-G(u) \cr
        &\geq \norm{\nabla u}_2^2 - N({\beta\over 2}-1)C_2(\norm u_\alpha^\alpha + \norm u_\beta^\beta). \cr}
    $$
By \ref[2.7]--\ref[2.8],
    $$  0 \geq \norm{\nabla u}_2^2 -N({\beta\over 2}-1)C_2C_3(\norm{\nabla u}_2^{\alpha'}
            + \norm{\nabla u}_2^{\beta'}).
    $$
Since $\alpha'$, $\beta'>2$, we have
    $$  \inf_{u\in \calP_m} \norm{\nabla u}_2^2 >0.         \eqno\label[2.9]
    $$
On the other hand, it follows from \ref[2.6] that $G(s)\leq {2\over \alpha-2}(\half g(s)s-G(s))$ for all $s\in\R$.
Thus for $u\in \calP_m$
    $$  \eqalignno{
        I(u) &\geq \half\norm{\nabla u}_2^2 -  {2\over \alpha-2}\intRN \half g(u)u -G(u) \cr
        &= \left( \half -{2\over (\alpha-2)N}\right) \norm{\nabla u}_2^2.  &\label[2.10] \cr}  
    $$
Since $\half -{2\over (\alpha-2)N}>0$ for $\alpha>2+{4\over N}$, \ref[2.9] and \ref[2.10] imply 
\ref[2.4].  \QED

\medskip

\proclaim Lemma \label[Lemma:2.2].
$I\in C^1(S_m,\R)$ satisfies $(PSP)_c$ for $c>0$.

\medskip

\claim Proof.
Suppose that $(u_j)_{j=1}^\infty\subset S_m$ satisfies
    $$  \eqalignno{
        &I(u_j)\to c,   &\label[2.11]\cr
        &\norm{dI(u_j)}_{(T_{u_j}S_m)^*} \to 0, &\label[2.12]\cr
        &P(u_j)\to 0.  &\label[2.13]\cr}
    $$

\Step 1:  $(u_j)$ is bounded in $\E$.

\smallskip

\noindent
Computing $\ref[2.11]-\half\cdot \ref[2.13]$, we have
    $$  -\left(1+{N\over 2}\right) \intRN G(u_j) +{N\over 4}\intRN g(u_j)u_j = c + o(1).
    $$
By (g2),
    $$  \left({N\over 4}\alpha-(1+{N\over 2})\right) \intRN G(u_j) \leq
         c+o(1) \leq
        \left({N\over 4}\beta-(1+{N\over 2})\right) \intRN G(u_j).   \eqno\label[2.14]
    $$
Thus, $\intRN G(u_j)$ and thus $\norm{\nabla u_j}_2^2$ are bounded.

\smallskip

\noindent
Since $(u_j)$ is bounded in $\E$, after taking a subsequence if necessary, we may assume
$u_j\wlimit u_0$ weakly in $\E$ and $u_j\to u_0$ strongly in $L^r(\R^N)$ for $r\in (2,2^*)$.
We note that
    $$  \intRN G(u_0)=\lim_{j\to\infty} \intRN G(u_j) >0                    \eqno\label[2.15]
    $$
follows from the second inequality in \ref[2.14].

\smallskip

\Step 2:  $(u_j)$ has a strongly convergent subsequence in $\E$.

\smallskip

\noindent
By \ref[2.12], there exists $(\lambda_j)\subset\R$ such that
    $$  -\Delta u_j +\lambda_j u_j -g(u_j) \to 0 \qquad \hbox{strongly in}\ (\E)^*. 
    $$
Thus,
    $$  \eqalign{
        \lambda_j m&= \lambda_j\norm{u_j}_2^2 \cr
        &=-\norm{\nabla u_j}_2^2 + \intRN g(u_j)u_j + o(1) \cr
        &= -P(u_j) - N\intRN \left(\half g(u_j)u_j -G(u_j)\right) + \intRN g(u_j)u_j + o(1)\cr
        &= \intRN \left( -{N-2\over 2} g(u_j)u_j +NG(u_j)\right) + o(1) \cr
        &\geq \left(-{N-2\over 2}\beta+N\right) \intRN G(u_j) +o(1). \cr}
    $$
We note $-{N-2\over 2}\beta+N>0$.  By \ref[2.15] we may assume $\lambda_j\to\lambda_0\in (0,\infty)$,
from which we deduce $u_j\to u_0$ strongly in $\E$.                          \QED

\medskip

\noindent
Now we apply our abstract deformation theory in Section \ref[Section:4] to our functional $I$.
We set $E=\E$ and $\Phi:\, \R\times E\to E$  by
    $$  (\Phi_\theta u)(x) = e^{{N\over 2}\theta} u(e^\theta x),
    $$
that is, $\Phi_\theta u= u_{e^\theta}$ using the notation \ref[2.2].
We also set
    $$  \eqalign{
        &S=S_m, \cr
        &J(\theta,u)=I(\Phi_\theta u) 
         = \half e^{2\theta}\norm{\nabla u}_2^2 - e^{-N\theta}\intRN G(e^{{N\over 2}\theta}u). \cr}
    $$
Note that the metric on $M=\R\times S$ is given by
    $$  \eqalign{
        \norm{(\kappa,v)}_{(\theta,u)} 
        &= (\abs\kappa^2 + \norm{\Phi_\theta v}_{H^1}^2)^{1/2} \cr
        &= (\abs\kappa^2 + e^{2\theta}\norm{\nabla v}_2^2 + \norm v_2^2)^{1/2} \cr
        &\qquad\hbox{for}\ (\theta,u)\in M=\R\times S \ \hbox{and for}\ (\kappa,v)\in T_{(\theta,u)}M.\cr}
    $$
It is easily observed that the assumption $(\Phi,S,I)$ is satisfied under these settings.  
We note that any solution $u$ of \ref[2.1], i.e., any critical point of $I(u)$, satisfies $P(u)=0$.
Thus by Proposition \ref[Proposition:4.5], we have for $K_c=\{ u\in S_m;\, I(u)=c,\, dI(u)=0 \ \hbox{on}\ 
T_u S_m\}$.

\medskip

\proclaim Proposition \label[Proposition:2.3].
For any $c>0$ and for any neighborhood $O$ of $K_c$ ($O=\emptyset$ if $K_c=\emptyset$) and 
any $\overline\epsilon>0$ there exist $\epsilon\in (0,\overline\epsilon)$ and $\eta\in 
C([0,1]\times S_m,S_m)$ such that (i)--(v) of Proposition \ref[Proposition:4.5] hold.  \QED

\medskip

\claim Remark \label[Remark:2.4].
In [\cite[BS1:7]], Bartsch and Soave take the natural constraint approach for \ref[2.1] under the
condition

\smallskip

\item{(g3)} A functional defined by $\wG(\xi)=\half g(\xi)\xi-G(\xi)$ is of class $C^1$ and satisfies
    $$  \wG'(\xi)\xi > (2+{4\over N})\wG(\xi) \quad \hbox{for}\ \xi\in\R\setminus\{ 0\}.
                                            \eqno\label[2.16]
    $$

\smallskip

\noindent
Under the condition (g3), they show that $\calP_m$ is a $C^1$-manifold and that $I$ restricted
to $\calP_m$ has properties corresponding to (P-i)--(P-ii) in Introduction.  These properties
enable them to show that for a suitable sequence of minimax methods for $I\big|_{\calP_m}$, 
the corresponding minimax values $c_n$ are actually critical values of $I$ on $S_m$ and satisfies $c_n\to \infty$
as $n\to\infty$.
However with their approach it seems difficult to obtain a multiplicity property 
as in Proposition 2.10 in Section \ref[Subsection:2.3] below.
In contrast, our Proposition \ref[Proposition:2.3] gives a deformation on $S_m$; we don't need to 
restrict $I$ on $\calP_m$, so we need not assume (g3). Moreover we can apply 
the genus theory (symmetric mountain pass theorem) to $I:\,S_m\to\R$.  
See Section \ref[Subsection:2.3].  

\medskip

\claim Remark \label[Remark:2.5].  
\ref[2.16] holds for $g(\xi)=\sum_{k=1}^\ell a_k\abs\xi^{p_k-2}\xi$, where $a_k>0$,
$p_k\in (2+{4\over N},2^*)$ ($k=1,2,\cdots,\ell$).

\medskip

\BSS{\label[Subsection:2.2]. Existence of a positive solution: Proof of (i) of Theorem \ref[Theorem:1.1]}
By Proposition \ref[Proposition:2.3], we can prove (i) of Theorem \ref[Theorem:1.1].

\medskip

\claim Proof of (i) of Theorem \ref[Theorem:1.1].
Applying the mountain pass theorem, we prove Theorem \ref[Theorem:1.1].

Let $b_0>0$ be given in Lemma \ref[Lemma:2.1] (ii).  We set
    $$  \eqalign{
    \Gamma =\{ \gamma\in C([0,1], S_m);\, &P(\gamma(0))<0<P(\gamma(1)), \cr
        &I(\gamma(0))<{b_0\over 2}, \ I(\gamma(1)) < {b_0\over 2}\}. \cr}
    $$
We note that $\Gamma\not=\emptyset$.  In fact, we fix $u_0\in S_m$ arbitrary and consider
$u_{0t}$ for $t\in (0,\infty)$.  By Lemma \ref[Lemma:2.1], we have for $L\gg 1$ and $0<\nu\ll 1$
    $$  \eqalign{
        &P(u_{0L}) < 0 < P(u_{0\nu}), \cr
        &I(u_{0L}) <{b_0\over 2}, \    I(u_{0\nu})< {b_0\over 2}. \cr}
    $$
Thus $\gamma(t)\equiv u_{0, \nu (1-\nu)+Lt} \in \Gamma$.   We set
    $$  b=\inf_{\gamma\in\Gamma}\max_{t\in [0,1]} I(\gamma(t)).     \eqno\label[2.17]
    $$
We can easily see that $\Gamma([0,1])\cap \calP_m \not=\emptyset$ for all $\gamma\in\Gamma$
and thus by Lemma \ref[Lemma:2.1] (ii),  we have $b\geq b_0$.

Thus, by Proposition \ref[Proposition:2.3], we can show that $b$ is a critical value of 
$I(u):\, S_m\to \R$ and \ref[2.1] has a solution.  \QED

\medskip

We note that by the proof of Lemma \ref[Lemma:2.2], we observe that if $u\in S_m$ is a critical point of
$I$, the corresponding Lagrange multiplier is positive.

\medskip

Next we shall prove $b=b_0$ when we suppose (g3) in addition to (g1)--(g2).  This result will be
important to study Schr\"odinger systems in Section \ref[Section:3].

\medskip

\proclaim Lemma \label[Lemma:2.6].
Assume that $N\geq 2$ and (g1)--(g3).
Then for $b>0$ given in \ref[2.17] and $b_0>0$ given in Lemma \ref[Lemma:2.1] (ii), we have $b=b_0$.

\medskip

\claim Proof.
By the proof of Theorem \ref[Theorem:1.1] (i), we know $b\geq b_0$.  We will show $b\leq b_0$ under (g3).

For an arbitrary fixed $u\in S_m$, we consider
    $$  I(u_t)
        =\half t^2 \norm{\nabla u}_2^2 - t^{-N}\intRN G(t^{N\over 2} u).
    $$
We have
    $$  \eqalign{
    {d\over dt }I(u_t) 
    &= t \norm{\nabla u}_2^2 - Nt^{-N-1}\intRN \wG(t^{N\over 2} u) \cr
    &= Nt\left( 
        {1\over N} \norm{\nabla u}_2^2 - t^{-N-2}\intRN\wG(t^{N\over 2}u) \right)\cr}
    $$
and we set $K(t)={1\over N}\norm{\nabla u}_2^2 - t^{-N-2}\intRN\wG(t^{N\over 2}u)$.
By (g3), we note that 
    $$  t\mapsto t^{-N-2}\wG(t^{N/2}\xi) \quad (\xi\not=0) 
    $$
is strictly increasing.
Therefore, if ${d\over dt}\Big|_{t=t_0}I(u_{t}) =0$, then $K(t_0)=0$ and thus we can see that $I(u_t)$ takes
a global maximum at $t=t_0$.

In particular, for any $u\in \calP_m$, we see $K(1)=0$,
    $$  \max_{t\in (0,\infty)} I(u_t)=I(u)
    $$
and
    $$  \gamma(t)\equiv u_{\nu(1-t)+Lt}\in \Gamma 
        \quad \hbox{for}\ L\gg 1 \ \hbox{and}\ 0<\nu\ll 1.
    $$
Thus we have $b\leq b_0$.  \QED

\medskip

\BSS{\label[Subsection:2.3].  Infinitely many solutions: Proof of (ii) of Theorem \ref[Theorem:1.1]}
We give a proof of (ii) of Theorem \ref[Theorem:1.1] using an idea related to symmetric mountain pass theorems
([\cite[R:26]], Chapter \raw[9]).

Recalling $u_t(x)=t^{N\over 2}u(tx)$, for $u\in\E\setminus\{ 0\}$ we define
    $$  h_0(u)=m^\half {u_{t(u)}\over \norm u_2}, \quad 
        \hbox{where}\ t(u)=\norm u_2.
    $$
We note that $h_0:\, \E\setminus\{ 0\}\to S_m$ is a continuous odd map.\m
We take a sequence $(E_k)_{k=1}^\infty$ of finite dimensional subspaces of $\E$ such that
    $$  \eqalign{
        &\dim E_k=k, \cr
        &E_1\subset E_2\subset\cdots\subset E_k\subset\cdots. \cr}
    $$
By (i) of Lemma \ref[Lemma:2.1] there exist sequences $(r_k)_{k=1}^\infty$, $(R_k)_{k=1}^\infty$ such that
    $$  \eqalign{
        &0<r_k<R_k, \cr
        &r_1>r_2>\cdots> r_k> \cdots >0, \cr
        &R_1<R_2<\cdots <R_k < \cdots, \cr}
    $$
and for $b_0>0$ given in Lemma \ref[Lemma:2.1]
    $$  I(u_{r_k})<{b_0\over 2}, \quad P(u_{r_k}) >0, \quad I(u_{R_k})<0, \quad P(u_{R_k})<0 \qquad
        \hbox{for all}\ u\in E_k\cap S_m.
    $$
We set
    $$  \eqalign{
        &D_k=\{ u\in E_k;\, r_k\leq \norm u_2\leq R_k\}, \cr
        &G_k =\{ h\in C(D_k, S_m);\, h(-u)=-h(u) \ \hbox{for all}\ u\in D_k, \cr
        &\qquad\qquad h(u)=h_0(u) \ \hbox{if}\ \norm u_2\in \{ r_k,R_k\}\}, \cr
        &\Gamma_j = \{ h(\overline{D_k\setminus Y});\, h\in G_k,\, k\geq j,\, Y\in \calE,\,
            \g(Y)\leq k-j\}. \cr}
    $$
Here $\calE$ is the family of sets $A\subset\E\setminus\{ 0\}$ such that $A$ is closed and symmetric
with respect to $0$.
For $A\in \calE$ the 
$\g(A)$ is defined as the smallest integer $n$ such that there exists
a continuous odd map $\varphi\in C(A,\R^n\setminus\{ 0\})$.  If there does not exist a finite such
$n$, we set $\g(A)=\infty$.  When $A=\emptyset$, we set $\g(\emptyset)=0$.

By our choice of $r_k$, $R_k$ and the definition of $h_0$,
    $$  \eqalign{
        &I(h_0(u))<{b_0\over 2}, \quad P(h_0(u))>0 \qquad \hbox{for}\ u\in D_k\ \hbox{with}\ \norm u_2=r_k,\cr
        &I(h_0(u))<0, \quad P(h_0(u))<0 \qquad \hbox{for}\ u\in D_k\ \hbox{with}\ \norm u_2=R_k.\cr}
    $$
Modifying the arguments in [\cite[R:26]], we have

\medskip

\proclaim Proposition \label[Proposition:2.7] {\rm (c.f. Proposition \raw[9.18] of [\cite[R:26]])}.
The sets $\Gamma_j$ have the following properties:
\item{(i)} $\Gamma_j\not=\emptyset$.
\item{(ii)} $\Gamma_{j+1}\subset\Gamma_j$.
\item{(iii)} If $\varphi\in C(S_m,S_m)$ is odd and $\varphi=id$ on $h_0(\partial D_k)$ for
all $k\geq j$.  Then $\varphi(B)\in\Gamma_j$ for all $B\in\Gamma_j$.
\item{(iv)} If $B\in\Gamma_j$, $Z\in\calE$ and $\g(Z)\leq s<j$, then $\overline{B\setminus Z}
\in \Gamma_{j-s}$. \QED

\medskip

\noindent
The following proposition gives an intersection property of $\Gamma_j$.

\medskip

\proclaim Proposition \label[Proposition:2.8] {\rm (c.f. Proposition \raw[9.23] of [\cite[R:26]])}.
For $j\in\N$, $B\in\Gamma_j$,
    $$  B\cap \calP_m\not=\emptyset.
    $$

\medskip

\claim Proof.
Set $B=h(\overline{D_k\setminus Y})$, where $k\geq j$, $h\in G_k$,  $\g(Y)\leq k-j$.
By our choice of $r_k$, $R_k$, we have for $u\in D_k$
    $$  \eqalign{
        &P(h(u))=P(h_0(u))>0 \quad \hbox{if}\ \norm u_2=r_k,\cr
        &P(h(u))=P(h_0(u))<0 \quad \hbox{if}\ \norm u_2=R_k.\cr}
    $$
Let $O$ be the connected
component of $\{ u\in D_k;\, P(h(u))>0\}$ including $\{u\in D_k;\, \norm u_2 = r_k\}$.
We note that $\wO\equiv O\cup \{ u\in D_k;\, \norm u_2\leq r_k\}$ is a bounded symmetric neighborhood of
$0$ in $E_k$.  Thus
    $$  \g(\partial \wO) =k.
    $$
It is easy to see that
    $$  h(\partial \wO) \subset \{ u\in S_m;\, P(u)=0\}. 
    $$
Set
    $$  W=\{ u\in D_k;\, P(h(u))=0\}.
    $$
We have $\partial \wO\subset W$ and $\g(W)=k$.  Thus
    $$  \g(\overline{W\setminus Y}) \geq k-(k-j)=j\geq 1.
    $$
In particular, $\overline{W\setminus Y}\not=\emptyset$.\m
On the other hand, we have  $h(\overline{W\setminus Y})\subset B\cap \calP$.  Thus
$B\cap \calP\not=\emptyset$.  \QED

\medskip

\noindent
Now we define
    $$  c_j=\inf_{B\in\Gamma_j}\max_{u\in B} I(u).
    $$
We have

\proclaim Corollary \label[Corollary:2.9].
\item{(i)} $c_1\leq c_2\leq \cdots\leq c_j\leq \cdots$.
\item{(ii)}  $c_j \geq b_0>0$ for all $j\in \N$, where $b_0$ is given in Lemma \ref[Lemma:2.1].

\medskip

\claim Proof.
(ii) follows from Proposition \ref[Proposition:2.8].  \QED

\medskip

\noindent
For $K_c$, we have

\medskip

\proclaim Proposition \label[Proposition:2.10] {\rm (c.f. Proposition \raw[9.30] of [\cite[R:26]])}.
If $c_j=c_{j+1}=\cdots=c_{j+p}\equiv d$, then
    $$  \g(K_d) \geq p+1.
    $$

\medskip

\claim Proof.  Since $I(u)$ satisfies $(PSP)_c$ for $c>0$, using our new deformation theory,
we can show Proposition \ref[Proposition:2.10] as in [\cite[R:26]]. \QED

\medskip

\proclaim Proposition \label[Proposition:2.11] {\rm (c.f. Proposition \raw[9.33] of [\cite[R:26]])}.
$c_j\to\infty$ as $j\to\infty$.  

\medskip

\claim Proof.
Following the argument for Proposition \raw[9.33] of [\cite[R:26]], we can show Proposition \ref[Proposition:2.11].
\QED

\medskip

\claim End of proof of (ii) of Theorem \ref[Theorem:1.1].
(ii) of Theorem \ref[Theorem:1.1] follows from Propositions \ref[Proposition:2.10] and \ref[Proposition:2.11].  \QED

\BScap{\label[Section:3]. Nonlinear Schr\"odinger systems}{\ref[Section:3] NS systems}
In this section we give another proof of Bartsch and Soave's result Theorem \ref[Theorem:1.2]
on nonlinear Schr\"odinger systems using deformation flow on $S_{m_1}\times S_{m_2}$.  

Since we consider the existence of positive solutions, setting $u_+=\max\{ u,0\}$, we study
the following system:
    $$  \left\{ \eqalign{
        -&\Delta u_1 +\lambda_1 u_1 = \mu_1 u_{1+}^3 + \beta u_1 u_2^2 \quad \hbox{in}\ \R^3, \cr
        -&\Delta u_2 +\lambda_2 u_2 = \mu_2 u_{2+}^3 + \beta u_1^2 u_2 \quad \hbox{in}\ \R^3, \cr
        &\norm{u_1}_2^2 = m_1, \ \norm{u_2}_2^2 = m_2, \cr}
        \right.                             \eqno\label[3.1]
    $$
where $\mu_i>0$, $m_i>0$ ($i=1,2$), $\beta<0$ are given constants and $\lambda_i\in\R$ ($i=1,2$) are
unknown Lagrange multipliers.  

To find a solution of \ref[3.1], we take a variational approach; we set
    $$  \eqalign{
        &S_{m_i}=\{ u\in H^1_r(\R^3);\, \norm u_2^2 =m_i\} \ (i=1,2), \cr
        &I_*(u_1,u_2) = \half\norm{\nabla u_1}_2^2 + \half \norm{\nabla u_2}_2^2
            -\int_{\R^3} G(u_1,u_2) :\, S_{m_1}\times S_{m_2}\to \R. \cr}
    $$
Here
    $$  G(u_1,u_2)={\mu_1\over 4}u_{1+}^4 + {\mu_2\over 4}u_{2+}^4 +{\beta\over 2}u_1^2u_2^2.
    $$
The Pohozaev functional for $I_*(u_1,u_2)$ is given by
    $$  P_*(u_1,u_2) = \norm{\nabla u_1}_2^2 + \norm{\nabla u_2}_2^2 - 3\int_{\R^3} G(u_1,u_2)
        :S_{m_1}\times S_{m_2}\to \R. 
    $$
We note that
    $$  P_*(u_1,u_2)={d\over dt}\Big|_{t=1} I(t^{3/2}u_1(tx), t^{3/2}u_2(tx)).
    $$
First we have

\medskip

\proclaim Proposition \label[Proposition:3.1].
Suppose that $(u_1,u_2)\in \SS$ is a critical point of $I_*$.
Then $u_1(x)>0$, $u_2(x)>0$ in $\R^3$ and for some $\lambda_1$, $\lambda_2>0$, \ref[3.1] holds.
Moreover the Pohozaev identity $P_*(u_1,u_2)=0$ holds.

\medskip

In next subsection, we prove Proposition \ref[Proposition:3.1] via a Liouville type argument.
As stated in Introduction we find a critical point $(u_1,u_2)$ of $I_*$ with the property
$P_*(u_1,u_2)=0$ via our new deformation argument for $I_*:\, \SS\to \R$ on $\SS$.

\medskip

\BSS{\label[Subsection:3.1].  Liouville type argument}
Here we develop a Liouville type argument for \ref[3.1] to prove Proposition \ref[Proposition:3.1].  Liouville
type argument is also important to verify the $(PSP)_c$ condition for $I_*(u_1,u_2)$.
See Proposition \ref[Proposition:3.8].
We remark that a similar result for positive solutions for \ref[1.2] is given in [\cite[BS1:7]].

\medskip

\proclaim Proposition \label[Proposition:3.2].
Suppose that $(u_1,u_2)\in \E\times \E$ and $\lambda_1$, $\lambda_2\in\R$ satisfy
$u_1\not=0$, $u_2\not=0$,
    $$  \eqalign{
        -&\Delta u_1 +\lambda_1 u_1 = \mu_1 u_{1+}^3 + \beta u_1 u_2^2 \quad \hbox{in}\ \R^3, \cr
        -&\Delta u_2 +\lambda_2 u_2 = \mu_2 u_{2+}^3 + \beta u_1^2 u_2 \quad \hbox{in}\ \R^3. \cr}
                                        \eqno\label[3.2]
    $$
Then $\lambda_1$, $\lambda_2>0$ and $u_1(x)$, $u_2(x)>0$ in $\R^3$.  Moreover $P_*(u_1,u_2)=0$ holds.

\medskip

\claim Proof.
First we note that $(u_1,u_2)$ satisfies 
    $$  \norm{\nabla u_1}_2^2 + \norm{\nabla u_2}_2^2 
        + \lambda_1\norm{u_1}_2^2 + \lambda_2\norm{u_2}_2^2 
            -4 \int_{\R^3} G(u_1,u_2) =0.       \eqno\label[3.3]
    $$
Next we remark $(u_1,u_2)$ also satisfies the Pohozaev identity:
    $$  \half\norm{\nabla u_1}_2^2 + \half\norm{\nabla u_2}_2^2 
        + 3\left({\lambda_1\over 2}\norm{u_1}_2^2 +{\lambda_2\over 2}\norm{u_2}_2^2 
            -\int_{\R^3} G(u_1,u_2) \right) =0.     \eqno\label[3.4]
    $$
In fact, by the argument in [\cite[BeL1:9], Proposition \raw[1]] we can show \ref[3.4].
We also note that $P_*(u_1,u_2)=0$ follows from \ref[3.3] and \ref[3.4].

Now we set $c=I_*(u_1,u_2)$ and we show $c>0$ and at least one of $\lambda_1$, $\lambda_2$ is positive.

Since $\lambda_1$, $\lambda_2$, $u_1$, $u_2$ satisfy \ref[3.3], $I_*(u_1,u_2)=c$ and $P_*(u_1,u_2)=0$,
we have
    $$  \eqalignno{
        &\norm{\nabla u_1}_2^2 + \norm{\nabla u_2}_2^2 =6c,     &\label[3.5]\cr
        &\int_{\R^3} G(u_1,u_2) = 2c,                               &\label[3.6]\cr
        &\lambda_1\norm{u_1}_2^2 + \lambda_2\norm{u_2}_2^2 = 2c.    &\label[3.7]\cr}
    $$
Since $u_1\not=0$, $u_2\not=0$, we have $c>0$ from \ref[3.5].  Thus by \ref[3.7], at least
one of $\lambda_1$, $\lambda_2$ must be positive.

In what follows, we assume that $\lambda_1>0$.  Then $u_1(x)$ is positive in $\R^3$ and decays
exponentially as $\abs x\to\infty$, that is, for some $c_1$, $c_2>0$
    $$  \abs{u_1(x)} \leq c_1 e^{-c_2\abs x} \quad \hbox{for all}\ x\in\R^3.
    $$
In fact, rewriting the first equation of \ref[3.2] as
    $$  -\Delta u_1 +(\lambda_1-\beta u_2^2) u_1 = \mu_1 u_{1+}^3 \quad \hbox{in}\ \R^3.
    $$
Noting $\beta<0$, we have the positivity and the decay property of $u_1(x)$.

If $\lambda_2>0$, in a similar way we can show that $u_2$ is positive and the conclusion of
Proposition \ref[Proposition:3.2] follows.  
Applying Proposition \ref[Proposition:3.3] to $\psi(x)=u_1(x)$ and $v(x)=u_2(x)$, we get
$\lambda_2>0$.  \QED

\medskip

\proclaim Proposition \label[Proposition:3.3].
Suppose that $\psi(x)\in H_r^1(\R^3)$ satisfies
    $$  \psi(x) = o({1\over \abs x}) \quad \hbox{as}\ \abs x\sim \infty.  \eqno\label[3.8]
    $$
For $\mu>0$, $\beta<0$, we consider
    $$  -\Delta v + (\lambda -\beta\psi^2)v = \mu v_+^3 \quad \hbox{in}\ \R^3.
                                                \eqno\label[3.9]
    $$
If \ref[3.9] has a non-zero solution $v\in H_r^1(\R^3)$, then $\lambda>0$.

\medskip

\claim Proof.
Suppose that $v(x)\in H_r^1(\R^3)$ satisfies \ref[3.9] and we show that 
$\lambda\leq 0$ cannot take a place.
We consider cases $\lambda=0$ and $\lambda<0$ separately.
Writing $r=\abs x$, we regard $\psi$, $v$ are functions of $r$.

\smallskip

{\sl
\Step 1:  Assume $\lambda=0$.  Then $v$ has finitely many zeros.

}

\smallskip

\noindent
We argue indirectly and assume that there exist $0<r_1<r_2<\cdots<r_n<r_{n+1}<\cdots$ such that
    $$  v(r_n)=0 \quad \hbox{and}\quad r_n\to \infty.
    $$
Setting $A_i=\{x\in\R^3;\, r_i<\abs x<r_{i+1}\}$, we have from \ref[3.9] that
    $$  \norm{\nabla v}_{L^2(A_i)}^2 = \mu\norm{v_+}_{L^4(A_i)}^4 + \beta\int_{A_i}\psi^2 v^2
        \leq \mu\norm v_{L^4(A_i)}^4.
    $$
By the Gagliard-Nirenberg inequality, there exists a constant $C>0$ such that
    $$  \norm u_4^4 \leq C \norm{\nabla u}_2^3\norm u_2 \quad \hbox{for}\ u\in H^1(\R^3).
    $$
Thus for $i=1,2,\cdots$
    $$  \norm{\nabla v}_{L^2(A_i)}^2 \leq C\norm{\nabla v}_{L^2(A_i)}^3\norm v_{L^2(A_i)},
    $$
from which we have
    $$  \eqalign{
        1 &\leq C \norm{\nabla v}_{L^2(A_i)} \norm v_{L^2(A_i)} \cr
        &\leq {C\over 2}(\norm{\nabla v}_{L^2(A_i)}^2 + \norm v_{L^2(A_i)}^2) \cr
        & = {C\over 2} \norm{\nabla v}_{H^1(A_i)}^2 \quad \hbox{for all}\ i\in\N. \cr}
    $$
Thus
    $$  \norm v_{H^1(\R^3)}^2 \geq \sum_{i=1}^\infty \norm v_{H^1(A_i)}^2=\infty,
    $$
which contradicts with $v\in H^1(\R^3)$.  Therefore $v(r)$ has only finitely many zeros.

By Step 1, there exists $R_0>0$ such that
    $$  v(x)\not=0 \quad \hbox{for}\ \abs x\geq R_0. 
    $$

{\sl
\Step 2: $\lambda=0$ cannot take a place.

}

\smallskip

\noindent
Here we use an idea from [\cite[BS1:7]].  We consider $\varphi(r)=r^{-\alpha}$ for $\alpha\in (1,{3\over 2}]$.
By the property \ref[3.8], it is easy to verify for some $R_1>R_0$
    $$  -\Delta\varphi -\beta\psi^2 \varphi < 0 \quad \hbox{for}\ \abs x\geq R_1.
    $$
First consider the case $v(r)>0$ for $\abs x\geq R_0$.  Since $-\Delta v-\beta\psi^2v=v_+^3$ in
$\abs x \geq R_1$, we have for $\delta>0$ small, $w(x)=v(x)-\delta\varphi(x)$ satisfies
    $$  \eqalign{
        -&\Delta w -\beta\psi^2 w>0 \quad \hbox{for}\ \abs x>R_1, \cr
        &w(R_1) >0, \cr
        &w(r)\to 0 \quad \hbox{as}\ r\to \infty.\cr}        \eqno\label[3.10]
    $$
Thus by the maximal principle, we have $w(x)>0$ for $\abs x\geq R_1$.  In particular, we have
    $$  v(x) \geq \delta \varphi(x) \quad \hbox{for}\ \abs x\geq R_1.
    $$
Noting $\varphi(x)\not\in L^2(\abs x\geq R_1)$, we have $v\not\in H^1(\R^3)$.  This is a contradiction.

Second we consider the case $v(x)<0$ for $\abs x\geq R_0$.  Since $-\Delta v -\beta \psi^2 v=0$
in $\abs x\geq R_1$, in a similar way, for small $\delta>0$ we can see $w(x)=-v(x)-\delta\varphi(x)$
satisfies \ref[3.10].  Thus $-v(x)\geq \delta \varphi(x)$ for $\abs x\geq R_1$ and we get a contradiction
again.
Thus $\lambda=0$ cannot take a place.

\smallskip

{\sl
\Step 3: $\lambda<0$ cannot take a place.

}

\smallskip

\noindent
Here we use an idea from [\cite[MMP:23], Lemma \raw[2.5], Step 4].  We set $w(r)=r v(r)$ and write
$'={d\over dr}$.   It follows from \ref[3.1] that
    $$  -w''  -[-\lambda+\beta\psi^2+\mu v_+^2]w =0 \quad \hbox{in}\ \R^3.
    $$
We set $V(r)=-\lambda+\beta\psi^2+\mu v_+^2$.  We have
    $$  \eqalignno{
        &V(r)\to -\lambda>0 \quad \hbox{as}\ r\to \infty,   &\label[3.11]\cr
        &V'(r)\in L^1(1,\infty).                            &\label[3.12]\cr}
    $$
In fact, it follows from $\psi$, $v\in H_r^1(\R^3)$ that $\psi(r)$, $v(r)\to 0$ as $r\to\infty$. 
Thus \ref[3.11] holds.
We also have $r^2\psi(r)^2$, $r^2\psi'(r)^2$, $r^2 v(r)$, $r^2 v'(r)^2\in L^1(0,\infty)$, from
which we deduce that $(\psi^2)'= 2\psi \psi'$, $(v^2)'=2vv'\in L^1(1,\infty)$ and thus \ref[3.12]
follows.

Next we set
    $$  E(r)= \half w'(r)^2 +\half V(r) w(r)^2.         \eqno\label[3.13]
    $$
By \ref[3.11], there exist $R_2>R_1$ and $C_1$, $C_2>0$ such that
    $$  C_1(w'(r)^2+w(r)^2) \leq E(r) \leq C_2(w'(r)^2+w(r)^2) \quad \hbox{for}\ r\geq R_2.
                                                            \eqno\label[3.14]
    $$
Differentiating \ref[3.13], we have
    $$  E'(r) = w'' w' + V(r)w w' + \half V'(r)w^2 =\half V'(r)w^2.
    $$
By \ref[3.14], $E'(r) \geq -{1\over 2C_1} \abs{V'(r)} E(r)$.  Thus we have
    $$  E(r) \geq E(R_2) \exp\left(-{1\over 2C_1}\int_{R_2}^r \abs{V'(s)}\, ds\right)
                        \quad \hbox{for} \ r\geq R_2,
    $$
In particular, by \ref[3.12], $\inf_{r\geq R_2} E(r)>0$.  By \ref[3.14] there exists $A_0>0$ such that
    $$  w'(r)^2+w(r)^2  \geq A_0 \quad \hbox{for} \ r\geq R_2.
    $$
By the definition of $w(r)$ and \ref[3.14],
    $$  \eqalign{
        \norm v_{H^1(\abs x\geq R_2)}^2
        &= c \int_{R_2}^\infty \left( v'(r)^2 + v(r)^2\right) r^2 \, dr \cr
        &= c \int_{R_2}^\infty \left( ((r^{-1}w(r))')^2 + r^{-2}w(r)^2\right) r^2 \, dr \cr
        &= c \int_{R_2}^\infty (w'(r)-r^{-1}w(r))^2 + w(r)^2\, dr \cr
        &= c \int_{R_2}^\infty w'^2+r^{-2}w^2 -2r^{-1}w w'  + w^2\, dr \cr
        &\geq c \int_{R_2}^\infty w'^2+r^{-2}w^2 - \half w'^2 -\half (2r^{-1}w)^2  + w^2\, dr \cr
        &= c \int_{R_2}^\infty \half w'^2+(1-r^{-2})w^2 \, dr \cr
        &=\infty.\cr}
    $$
It is a contradiction and $\lambda<0$ cannot take a place.

Thus we complete the proof of Proposition \ref[Proposition:3.3].  \QED

\medskip

\claim Proof of Proposition \ref[Proposition:3.1].
Let $(u_1,u_2)\in \SS$ be a critical point of $I_*:\,\SS\to\R$.  It is clear that \ref[3.2] holds
for some $\lambda_1$, $\lambda_2\in \R$.  Hence, the desired result follows from
Proposition \ref[Proposition:3.2].  \QED

\medskip

\claim Remark \label[Remark:3.4].
Modifying the proof of Proposition \ref[Proposition:3.1], we can show that if $(u_1,u_2)\in (\E\setminus\{0\})^2$
and $\lambda_1$, $\lambda_2\in\R$ satisfy
    $$  \left\{ \eqalign{
        -&\Delta u_1 +\lambda_1 u_1 = \mu_1 u_{1}^3 + \beta u_1 u_2^2 \quad \hbox{in}\ \R^3, \cr
        -&\Delta u_2 +\lambda_2 u_2 = \mu_2 u_{2}^3 + \beta u_1^2 u_2 \quad \hbox{in}\ \R^3, \cr}
        \right.                             
    $$
then we have $\lambda_1$, $\lambda_2>0$.

\medskip


\BSS{\label[Subsection:3.2]. A minimax method for $I_*$}
We also define for $i=1,2$
    $$  \eqalign{
        I_i(u)= \half\norm{\nabla u}_2^2 -{\mu_i\over 4}\norm{u_+}_4^4:\, S_{m_i}\to \R, \cr
        P_i(u)= \norm{\nabla u}_2^2 -{3\mu_i\over 4}\norm{u_+}_4^4:\, S_{m_i}\to \R. \cr}
    $$
We have

\medskip

\proclaim Lemma \label[Lemma:3.5].
For $i=1,2$, $I_i:\, S_{m_i}\to \R$  has a unique critical point and 
the Lagrange multiplier $\lambda_i$ is positive.

\claim Proof.  Using the Pohozaev identity, we have $\lambda_i>0$.
We note that for any given $\lambda>0$,
    $$  -\Delta u + \lambda u = \mu_i u_+^3 \ \hbox{in}\ \R^3, \quad u\in H^1_r(\R^N)
    $$
has a unique solution $u(\lambda;x)=({\lambda\over\mu_i})^{1/2}\omega(\lambda^{1/2}x)$.  
Here $\omega$ is a unique solution of
    $$  -\Delta \omega + \omega = \omega_+^3 \ \hbox{in}\ \R^3, \quad \omega\in H^1_r(\R^N).
    $$  
We know that
    $$  \norm{\nabla\omega}_2^2= 3 \norm{\omega}_2^2, \quad \norm\omega_4^4=4 \norm\omega_2^2.
    $$
Thus
    $$  \norm{u(\lambda;x)}_2^2 = {1\over \mu_i\lambda^{1/2}}\norm\omega_2^2, \quad
        I_i(u(\lambda;x))=\half {\lambda^{1/2}\over \mu_i}\norm\omega_2^2.
    $$
Thus for given $m_i>0$ and $\mu_i>0$, there is a unique $\lambda>0$ such that 
$\norm{u(\lambda;x)}_2^2 = m_i$.  \QED

\medskip

We also denote the unique critical value of $I_i(u)$ by $b_i$.  By Theorem \ref[Theorem:1.1] (i) and
Lemma \ref[Lemma:2.1],
    $$  b_i= \inf\{ I_i(u);\, u\in S_{m_i}, \ P_i(u)=0\}>0.
    $$
By the assumption $\beta<0$, we have
    $$  I_*(u_1,u_2) \geq I_1(u_1) + I_2(u_2) \quad \hbox{for all}\ (u_1,u_2)\in S_{m_1}\times S_{m_2}.
    $$
We introduce the following minimax value:
    $$  \eqalign{
        b_* &=\inf_{\gamma\in \Gamma_*}\max_{(s,t)\in [0,1]^2} I_*(\gamma(s,t)), \cr
        \Gamma_* &= \{\gamma(s,t)=(\gamma_1(s,t), \gamma_2(s,t))
            \in C([0,1]^2, S_{m_1}\times S_{m_2});\, \cr
            &\qquad\quad P_1(\gamma_1(0,t))>0> P_1(\gamma_1(1,t)) \ \hbox{for}\ t\in [0,1], \cr
            &\qquad\quad   I_1(\gamma_1(0,t)), \  I_1(\gamma_1(1,t)) < b_1 \ \hbox{for}\ t\in [0,1], \cr
            &\qquad\quad P_2(\gamma_2(s,0))>0> P_2(\gamma_2(s,1)) \ \hbox{for}\ s\in [0,1], \cr
            &\qquad\quad   I_2(\gamma_2(s,0)), \  I_2(\gamma_2(s,1)) < b_2 \ \hbox{for}\ s\in [0,1], \cr
            &\qquad\quad I_*(\gamma(s,t)) < \overline b \ \hbox{for} \ (s,t)\in \partial([0,1]^2)\}, \cr}
    $$
where
    $$  \overline b \in (\max\{ b_1,b_2\}, b_1+b_2).
    $$
We note that

\medskip

\proclaim Lemma \label[Lemma:3.6].
$\Gamma_*\not=\emptyset$.

\medskip

\claim Proof.
For given $\delta>0$, we choose $w_i(x)\in S_{m_i}\cap \{P_i(u)=0\}\cap C_0^\infty(\R^3)$ ($i=1,2$)
such that
    $$  I_i(w_i)\in [b_i,b_i+\delta] \qquad \hbox{and}\qquad 0\not\in\supp w_i.
    $$
Setting $w_{it}(x)=t^{3/2}w_i(tx)$ for $t>0$, by Lemma \ref[Lemma:2.1], we have for $i=1,2$
    $$  \eqalign{
        &I_i(w_{it}),\, \norm{\nabla w_{it}}_2^2,\, \norm{w_{it}}_4^4\to +0 \quad
        \hbox{as}\ t\to +0,\cr
        &I_i(w_{it}) \to -\infty \quad \hbox{as}\ t\to \infty,\cr
        &I_i(w_i)=I_i(w_{i1})= \sup_{t\in (0,\infty)} I_i(w_{it}), \cr
        &P_i(w_{it}) > 0 \qquad \hbox{for}\ t\in (0,1),\cr
        &P_i(w_{it}) < 0 \qquad \hbox{for}\ t\in (1,\infty).\cr}
    $$
Note that
    $$  \int_{\R^3} w_{1s}^2 w_{2t}^2 
        \leq \cases{\norm{w_{1s}}_\infty^2 m_2 \cr \norm{w_{2t}}_\infty^2 m_1 \cr} 
        = \cases{ s^{3/2} \norm{w_{1}}_\infty^2 m_2 \cr t^{3/2} \norm{w_{2}}_\infty^2 m_1 \cr} 
    $$
and
    $$  \eqalign{
        I_*(w_{1s},w_{2t}) &\leq I_1(w_{1s}) + I_2(w_{2t}) + {\abs \beta\over 2}\int_{\R^3} w_{1s}^2 w_{2t}^2 \cr
        &\leq I_1(w_{1s}) + I_2(w_{2t}) + {\abs \beta\over 2}
        \cases{ s^{3/2} \norm{w_{1}}_\infty^2 m_2 \cr t^{3/2} \norm{w_{2}}_\infty^2 m_1 \cr}. \cr}
    $$
We choose $L_i>1$ and $\nu_i\in (0,1)$ such that
    $$  \eqalignno{
        &I_*(w_{1\nu_1},w_{2\nu_2}) <\delta, \cr
        &I_*(w_{1s},w_{2\nu_2}) < b_1+2\delta \quad \hbox{for}\ s\in [\nu_1,L_1],\cr
        &I_*(w_{1\nu_1},w_{2t}) < b_2+2\delta \quad \hbox{for}\ t\in [\nu_2,L_2],\cr
        &I_1(w_{1L_1})<0, \ I_2(w_{2L_2})<0, \cr
        &I_*(w_{1L_1},w_{2\nu_2})<0, \ I_*(w_{1\nu_1},w_{2L_2})<0, \ I_*(w_{1L_1},w_{2L_2})<0.  &\label[3.15]\cr}
    $$
Choosing $\nu_i$ smaller and $L_i>1$ larger if necessary, we may assume
    $$  \supp w_{1\nu_1},\ \supp w_{1L_1},\ \supp w_{2\nu_2},\ \supp w_{2L_2} \quad 
        \hbox{are pair-wise disjoint.}  \eqno\label[3.16]
    $$
Setting
    $$  \gamma_0(s,t)=(w_{1, \nu_1+(L_1-\nu_1)s}, w_{2, \nu_2+(L_2-\nu_2)t}) \quad
        \hbox{for}\ (s,t)\in (\{0\}\times [0,1])\cup ([0,1]\times \{ 0\}),
    $$
we observe that
$\gamma_0(s,t)$ possesses the desired properties on $(\{0\}\times [0,1])\cup ([0,1]\times \{ 0\})$
in the definition of $\Gamma_*$.

In what follows, we define $\gamma(s,t)$ on $(\{1\}\times [0,1])\cup ([0,1]\times \{ 1\})$ in 3 steps.

\smallskip

\Step 1:  There exists a continuous path $c(s):\, [0,1]\to S_{m_1}\times S_{m_2}$ joining
$(w_{1\nu_1},w_{2L_2})$, $(w_{1L_1},w_{2L_2})$ and
$\int_{\R^3} G(c(s))>0$ for all $s\in [0,1]$.

\smallskip

\noindent
In fact, setting
    $$  c(s)=\left( \sqrt{m_1} {(1-s)w_{1\nu_1}+sw_{1L_1} \over \norm{(1-s)w_{1\nu_1}+sw_{1L_1}}_2},
        w_{2L_2}\right).
    $$
it follows from \ref[3.16] that $\int_{\R^3} G(c(s)) >0$ for $s\in [0,1]$.

\smallskip

\Step 2:  There exists a continuous path $\tilde c(s)=(\tilde c_1(s), \tilde c_2(s)):\, [0,1]\to S_{m_1}\times S_{m_2}$ joining
$(w_{1\nu_1},w_{2L_2})$, $(w_{1L_1},w_{2L_2})$ and
    $$  \eqalignno{
        &I_*(\tilde c(s)) \leq 0 \qquad \hbox{for}\ s\in [0,1], &\label[3.17]\cr
        &P_2(\tilde c_2(s)) < 0 \qquad \hbox{for}\ s\in [0,1].  &\label[3.18]\cr}
    $$

\smallskip

\noindent
We note that for $c(s)=(c_1(s),c_2(s))$ obtained in Step 1,
    $$  \eqalign{
        &I_*(c(s)_t)= {t^2\over 2}(\norm{\nabla c_1(s)}_2^2+ \norm{\nabla c_2(s)}_2^2)
            -{t^3\over 4}\int_{\R^3} G(c(s)), \cr
        &P_2(c_2(s)_t)= {t^2}\norm{\nabla c_2(s)}_2^2 - {3\over 4}\mu_2 t^4 \norm{c_2(s)}_4^4.\cr}
    $$
Since $\int_{\R^3} G(c(s))>0$ for $s\in [0,1]$, we have $I_*(c(s)_t)<0$ for large $t>1$.  We set
    $$  \tilde c(s)=c(s)_{T(s)},
    $$
where 
    $$  T(s)=\inf\{ t\in [1,\infty);\, I_*(c(s)_t)<0\} 
        = \max\left\{ 1, {2(\norm{\nabla c_1(s)}_2^2 + \norm{\nabla c_2(s)}_2^2)\over \int_{\R^3} G(c(s))} \right\}.
    $$
By \ref[3.15] we observe that $\tilde c(0)=(w_{1\nu_1}, w_{2L_2})$, $\tilde c(1)=(w_{1L_1}, w_{2L_2})$ and
$\tilde c(s)$ has the desired properties \ref[3.17]--\ref[3.18].

\smallskip

\Step 3:  Conclusion.

\smallskip

\noindent
In a similar way to Steps 1--2, we can find a path $\hat c(t)=(\hat c_1(t),\hat c_2(t)):\, [0,1]\to S_{m_1}\times S_{m_2}$ joining
$(w_{1L_1},w_{2\nu_2})$, $(w_{1L_1},w_{2L_2})$ and
    $$  \eqalign{
        &I_*(\hat c(t)) \leq 0 \qquad \hbox{for}\ t\in [0,1], \cr
        &P_1(\hat c_1(t)) < 0 \qquad \hbox{for}\ t\in [0,1]. \cr}
    $$
We set $\gamma(s,t):\, \partial ([0,1]\times[0,1])\to S_{m_1}\times S_{m_2}$ by
    $$  \gamma(s,t)=\cases{
        \gamma_0(s,t) &for $(s,t)\in (\{0\}\times [0,1])\cup([0,1]\times \{0\})$, \cr
        \tilde c(s) &for $(s,t)\in [0,1]\times \{1\}$, \cr
        \hat c(t) &for $(s,t)\in \{1\}\times [0,1]$. \cr}
    $$
Extending $\gamma(s,t)$ continuously on $[0,1]\times [0,1]$, we see $\gamma\in\Gamma_*$ and
$\Gamma_*\not=\emptyset$. \QED

\medskip

\noindent
We have

\proclaim Lemma \label[Lemma:3.7].  
$b_* \geq b_1+b_2$.

\claim Proof.
Using the degree theory, for any $\gamma\in \Gamma_*$ we find $(t_{01},t_{02})\in [0,1]^2$ such
that
    $$  P_1(\gamma(t_{01},t_{02})) = P_2(\gamma(t_{01},t_{02})) = 0,
    $$
from which we have $b_* \geq b_1+b_2$.  \QED

\medskip

\BSS{\label[Subsection:3.3]. $(PSP)$ condition}
Next we show the following proposition, which is a key to generate our new deformation flow.

\proclaim Proposition \label[Proposition:3.8].  
For $c>\max\{ b_1,b_2\}$, $I(u_1,u_2)$ satisfies $(PSP)_c$ condition on $S_{m_1}\times S_{m_2}$.

\claim Proof.
Assume $(U_j)_{j=1}^\infty=(u_{1j},u_{2j})_{j=1}^\infty \subset S_{m_1}\times S_{m_2}$ satisfies
as $j\to\infty$
    $$  \eqalignno{
        &I_*(U_j)\to c > \max\{ b_1,b_2\}, &\label[3.19]\cr
        &\norm{I_*'(U_j)}_{(T_{U_j}(S_{m_1}\times S_{m_2}))^*} \to 0, &\label[3.20] \cr
        &P_*(U_j)\to 0, &\label[3.21]\cr}
    $$
We note that \ref[3.20] implies for some $(\lambda_{1j},\lambda_{2j})\in\R^2$
    $$  \eqalignno{
        &-\Delta u_{1j} +\lambda_{1j}u_{1j} = \mu_1 u_{1j+}^3 + \beta u_{1j}u_{2j}^2 + o(1)\norm{U_j}_{H^1}, &\label[3.22]\cr
        &-\Delta u_{2j} +\lambda_{2j}u_{2j} = \mu_2 u_{2j+}^3 + \beta u_{1j}^2u_{2j} + o(1)\norm{U_j}_{H^1}. &\label[3.23]\cr}
    $$

\smallskip

\Step 1:  $(U_j)_{j=1}^\infty$ is bounded in $H^1_r(\R^3)\times H_r^1(\R^3)$.

\smallskip

\noindent
It follows from \ref[3.19] and \ref[3.21] that
    $$  I_*(U_j)= {1\over 6}(\norm{\nabla u_{1j}}_2^2 + \norm{\nabla u_{2j}}_2^2) + o(1) = c+o(1).
    $$
Thus we have boundedness of $(U_j)_{j=1}^\infty$ in $H^1(\R^3)\times H^1(\R^3)$.

\medskip

After taking a subsequence if necessary,  we may assume that $u_{1j}\wlimit u_{10}$, $u_{2j}\wlimit u_{20}$
weakly in $H^1_r(\R^3)$.  We note that $u_{1j}\to u_{10}$, $u_{2j}\to u_{20}$ strongly in
$L^4(\R^3)$.

\medskip

\Step 2:  $(\lambda_{1j}, \lambda_{2j})$ is bounded in $\R^2$ and we may assume $\lambda_{1j}\to\lambda_{10}$
and $\lambda_{2j}\to\lambda_{20}$.  Moreover we have
    $$  \lambda_{10}m_1 + \lambda_{20}m_2 = 2c.         \eqno\label[3.24]
    $$
In particular, at least one of $\lambda_{10}$, $\lambda_{20}$ is positive.

\smallskip

\noindent
Multiplying $u_{1j}$ to \ref[3.22] and multiplying $u_{2j}$ to \ref[3.23], we can easily see boundedness of
$(\lambda_{1j}, \lambda_{2j})$.  By \ref[3.22] and \ref[3.23],
    $$  \norm{\nabla u_{1j}}_2^2 + \norm{\nabla u_{2j}}_2^2 -4\int_{\R^3} G(u_{1j},u_{2j})
        + \lambda_{1j}m_1 + \lambda_{2j}m_2 =o(1).
    $$
By \ref[3.21], we have
    $$ \lambda_{1j}m_1 + \lambda_{2j}m_2 = {1\over 3} (\norm{\nabla u_{1j}}_2^2 + \norm{\nabla u_{2j}}_2^2)
        = 2c +o(1).
    $$
Thus \ref[3.24] follows.

\medskip

\Step 3:  $u_{10}\not=0$ and $u_{20}\not=0$.

\smallskip

\noindent
Assume that $u_{20}=0$.  By \ref[3.19] and \ref[3.21], we have
    $$  I_*(U_j) = \half \int_{\R^3} G(u_{1j}, u_{2j}) \to {\mu_1\over 8}\norm{u_{10+}}_4^4.
    $$
In particular,  we have
    $$  {\mu_1\over 8}\norm{u_{10+}}_4^4 = c        \eqno\label[3.25]
    $$
and $u_{10}$ is non-trivial.  By \ref[3.23], we also have
    $$  \norm{\nabla u_{2j}}_2^2 + \lambda_{2j} m_2 \to 0,
    $$
from which we deduce $\lambda_{20}\leq 0$ and thus $\lambda_{10}>0$.

By \ref[3.22], we also have
    $$  \norm{\nabla u_{1j}}_2^2 + \lambda_{10}\norm{u_{1j}}_2^2 = \mu_1\norm{u_{1j}}_4^4 +o(1).
    $$
Since $\lambda_{10}>0$, we easily deduce that $u_{1j}\to u_{10}$ strongly in $H_r^1(\R^3)$.
Thus $u_{10}\in S_{m_1}$ is a critical point of $I_1(u)$ under constraint $\norm u_2^2=m_1$.
By the uniqueness of critical points of $I_1(u)$, we have
    $$  b_1=I_1(u_{10})= {\mu_1\over 8}\norm{u_{10+}}_4^4.
    $$
Since $c> b_1$, it contradicts \ref[3.25] and we have $u_{20}\not=0$.  In a similar way we can show
$u_{10}\not=0$.
    
\medskip

\Step 4: $\lambda_{10}>0$, $\lambda_{20}>0$ and $(u_{10},u_{20})\in\SS$.  Moreover $u_{1j}\to u_{10}$,
$u_{2j}\to u_{20}$ strongly in $H_r^1(\R^3)$.  We also have $dI_*(u_{10},u_{20})=0$ and $P_*(u_{10},u_{20})=0$.

\smallskip

\noindent
By \ref[3.24], at least one of $\lambda_{10}$, $\lambda_{20}$ is positive.  We assume $\lambda_{10}>0$.
Then we can see that $u_{10}>0$ in $\R^3$ and it decays exponentially as $\abs x\to\infty$.
Applying Proposition \ref[Proposition:3.3] to 
    $$  -\Delta u_{20} +\lambda_{20}u_{20}=\mu_2 (u_{20+})^3 +\beta u_{10}^2 u_{20} \quad \hbox{in}\ \R^3,
    $$
we have $\lambda_{20}>0$.

Since $\lambda_{10}$, $\lambda_{20}>0$, it is not difficult to see that $u_{1j}\to u_{10}$,
$u_{2j}\to u_{20}$ strongly in $\E$, from which we also deduce that $dI_*(u_{10},u_{20})=0$ and 
$P_*(u_{10},u_{20})=0$.  Thus $(PSP)_c$ holds.  \QED

\medskip

Our Proposition \ref[Proposition:3.8] enables us to apply our abstract deformation theory to $I_*(u_1,u_2)$.
More precisely, we set $E=H_r^1(\R^3)\times H_r^1(\R^3)$ and define
$\Phi:\, \R\times E \to E$ by
    $$  \Phi_\theta (u_1,u_2)(x)
        = (e^{{3\over 2}\theta} u_1(e^\theta x), e^{{3\over 2}\theta} u_2(e^\theta x)).
    $$
We also set
    $$  \eqalign{
        &S=\SS, \cr
        &J(\theta,u_1,u_2) = I_*(\Phi_\theta(u_1,u_2)) 
            =\half e^{2\theta}(\norm{\nabla u_1}_2^2+\norm{\nabla u_2}_2^2)
                - e^{3\theta} \int_{\R^3} G(u_1,u_2). \cr}
    $$
It is easily observed that the assumption $(\Phi,S,I)$ is satisfied under these settings.  Thus
by Proposition \ref[Proposition:4.5], we have Proposition \ref[Proposition:1.3] in Introduction.  Since $(PSP)_c$ holds for
$c>\max\{ b_1,b_2\}$, we have

\medskip

\proclaim Proposition \label[Proposition:3.9].
For any $c>\max\{ b_1,b_2\}$ and for any neighborhood $O$ of $K_c$ ($O=\emptyset$ if $K_c=\emptyset$) and 
any $\overline\epsilon>0$ there exist $\epsilon\in (0,\overline\epsilon)$ and $\eta\in 
C([0,1]\times \SS, \SS)$ such that (i)--(v) of Proposition \ref[Proposition:1.3] hold.  \QED

\medskip

\noindent
Now we can complete the proof of Theorem \ref[Theorem:1.2].  

\medskip

\claim End of the proof of Theorem \ref[Theorem:1.2].  
By Proposition \ref[Proposition:3.9], our new deformation result on $S_{m_1}\times S_{m_2}$ enables us
to show that $b_* (\geq b_1+b_2)$ is a critical value of $I_*(U)$.  
Thus \ref[3.1] has at least one solution $(u_1, u_2)$.  
\QED

\medskip

\BS{\label[Section:4]. Deformation argument}
To give proofs of Propositions \ref[Proposition:1.3] and \ref[Proposition:2.3] systematically, we give our deformation result in an abstract
setting.

Let $(E,\norm\cdot_E)$ be a Banach space and let $\Phi_\theta:\, \R \to L(E); \theta\mapsto \Phi_\theta$ 
be a continuous group action of $\R$ and  we suppose there exists an embedded $C^2$-submanifold $S$ of
$E$ and $I(u)\in C^1(S,\R)$ which satisfy the following assumptions:

\medskip

\claim Assumption $(\Phi,S,I)$.
\item{(i)} $\Phi_\theta$ is $C^0$ group action, that is, 
    $$  \eqalign{
        &\Phi_\theta\circ \Phi_{\theta'} = \Phi_{\theta+\theta'} \quad\hbox{for}\ \theta, \theta'\in \R, \cr
        &\Phi_0=id, \cr 
        &\theta\mapsto \Phi_\theta u \ \hbox{is strongly continuous for any}\ u\in E.\cr}
    $$
\item{(ii)} $S$ is invariant under $\Phi_\theta$, that is, $\Phi_\theta(S)\subset S$ for all $\theta\in\R$.
\item{(iii)}  Let $M=\R\times S$ and on the tangent bundle 
    $$  TM=\coprod_{(\theta,u)\in M} (\R\times T_uS),
    $$
we introduce a metric
    $$  \norm{(\kappa,v)}_{(\theta,u)} = \left(\kappa^2+\norm{\Phi_\theta v}_E^2\right)^{1/2}.  
                                                                            \eqno\label[4.1]
    $$
We assume $\norm\cdot_{(\theta,u)}$ is a metric of class $C^2$ on $TM$.
\item{(iv)} Let 
    $$  J(\theta,u)=I(\Phi_\theta u):\, M=\R\times S\to \R.     \eqno\label[4.2]
    $$
We assume that $J(\theta,u)$ is of class $C^1$ on $M$.

\medskip

\noindent
We remark Assumption $(\Phi,S,I)$ holds in rather special settings.
We give examples, which cover \ref[1.1], \ref[1.2] and results in [\cite[HT:16]].

\medskip

\claim Example \label[Example:4.1].
In the setting of Sections \ref[Section:2]--\ref[Section:3], Assumption $(\Phi,S,I)$ holds.

\medskip

\claim Example \label[Example:4.2] ([\cite[HT:16]]).
Let $E=H_r^1(\R^N)$ ($N\geq 2$) and let
    $$  (\Phi_\theta u)(x) = u(x/e^\theta).
    $$
Then $\Phi_\theta$ is a $C^0$ group action of $\R$ (not of class $C^1$).  Set $S=H_r^1(\R^N)$.
Then $M=\R\times H_r^1(\R^N)$ and 
    $$  \eqalign{
    \norm{(\kappa,v)}_{(\theta,u)}^2 &=\kappa^2 + \norm{\Phi_\theta v}_{H^1(\R^N)}^2 \cr
    &= \kappa^2 + e^{(N-2)\theta}\norm{\nabla v}_{L^2(\R^N)}^2 + e^{N\theta}\norm v_{L^2(\R^N)}^2. \cr}
    $$
gives a $C^2$ metric.  Under conditions (g1)--(g2),
we consider
    $$  I(u) = \half\norm{\nabla u}_{L^2(\R^N)}^2 -\int_{\R^N}G(u):\, H_r^1(\R^N)\to \R.
    $$
Then
    $$  J(\theta,u) = I(\Phi_\theta u) 
        = {e^{(N-2)\theta}\over 2}\norm{\nabla u}_{L^2(\R^N)}^2 -e^{N\theta}\int_{\R^N}G(u)
    $$
is of class $C^1$.

\medskip

\noindent
Under the assumption $(\Phi,S,I)$, we set
    $$  P(u)=\partial_\theta J(0,u):\, S\to \R,
    $$
which corresponds to the Pohozaev functional.

We denote by $dI(u)$ the derivative of $I$ and by $\norm{dI(u)}_{T_u^*S}$ its norm, that is,
    $$  \norm{dI(u)}_{T_u^*S} = \sup_{v\in T_uS, \norm v_{T_uS}\leq 1} dI(u)v.  \eqno\label[4.3]
    $$
We also impose the following Palais-Smale type condition.

\medskip

\claim Definition \label[Definition:4.3].
For $b\in\R$, we say that $I(u)$ satisfies $(PSP)_b$ condition on $S$ if any sequence
$(u_n)_{n=1}^\infty\subset S$ with
    $$  \eqalign{
        &I(u_n)\to b,                       \cr
        &\norm{dI(u_n)}_{T_{u_n}^*S}\to 0,  \cr
        &P(u_n)\to 0,                       \cr}
    $$
has a strongly convergent subsequence.

\medskip

In what follows, we use the following notation: for $c\in\R$
    $$  \eqalign{
        &[I\leq c]_S=\{ u\in S;\, I(u)\leq c\}, \cr
        &\widehat K_c = \{ u\in S;\, I(u)=c,\, dI(u)=0,\, P(u)=0\}. \cr}
    $$
We note that the definition of our critical set $\widehat K_c$ is different from the standard one, that is,
we require $P(u)=0$ in addition to $dI(u)=0$.

\medskip

\claim Remark \label[Remark:4.4].
In the settings of Sections \ref[Section:2]--\ref[Section:3], any critical point satisfies the Pohozaev identity $P(u)=0$.  That is,
	$$	\widehat K_c=K_c
	$$
holds, where $K_c=\{ u\in S;\, I(u)=c,\, dI(u)=0,\}$.

\medskip

The aim of this section is to show the following deformation result.

\medskip

\proclaim Proposition \label[Proposition:4.5].
Suppose that assumption $(\Phi,S,I)$ and $(PSP)_b$ for $b\in\R$ hold.
For any neighborhood $O$ of $\widehat K_b$ ($O=\emptyset$ if $\widehat K_b=\emptyset$) and any $\overline\epsilon>0$,
there exist $\epsilon\in (0,\overline\epsilon)$ and $\eta\in C([0,1]\times S,S)$ such that
\item{(i)} $\eta(0,u)=u$ for $u\in S$.
\item{(ii)} $\eta(t,u)=u$ for $t\in [0,1]$ if $u\in [I\leq b-\overline\epsilon]_S$.
\item{(iii)} $t\mapsto I(\eta(t,u))$ is non-increasing for $u\in S$.
\item{(iv)} $\eta(1,[I\leq b+\epsilon]_S\setminus O) \subset [I\leq b-\epsilon]_S$, 
$\eta(1,[I\leq b+\epsilon]_S) \subset [I\leq b-\epsilon]_S\cup O$.

{\sl \noindent
Moreover, if $S$ is symmetric with respect to $0$ and $I(u)$ is even in $u$, that is,
    $$  \eqalignno{
        &-S=S, \quad 0\not\in S,                    &\label[4.4] \cr
        &I(-u)=I(u) \quad \hbox{for}\ u\in S,       &\label[4.5]\cr}
    $$
then we also have 
\item{(v)} $\eta(t,-u)=-\eta(t,u) \quad \hbox{for}\ (t,u)\in [0,1]\times S$.
}

\medskip

\claim Remark \label[Remark:4.6].
Under the assumption $(\Phi,S,I)$, if a value $b$ is given by a minimax method and $(PSP)_b$ holds,
Proposition \ref[Proposition:4.5] implies that $\widehat K_b\not=\emptyset$.  That is, there exists a critical point
$u$ of $I$ with the property $P(u)=0$.

In general, for example for nonlinear equations involving fractional operators, it is difficult to
check $\widehat K_c=K_c$.  In such a situation our Proposition \ref[Proposition:4.5] ensures the existence of a 
critical point with the Pohozaev property $P(u)=0$.

As another advantage of our approach, Proposition \ref[Proposition:4.5] can be applied to obtain multiplicity result
as in Section \ref[Subsection:2.3].  We note that approaches in [\cite[BdV:6]] and [\cite[HIT:17]] ensure just the existence
of a Palais-Smale sequence at some minimax levels and it seems difficult to use benefits of topological
tools like the genus (e.g. Proposition \ref[Proposition:2.10]) directly.  As we show in Section \ref[Subsection:2.3], our deformation
result works well together with the genus theory.

\medskip

To show Proposition \ref[Proposition:4.5], as in [\cite[HT:16]], we exploit the functional $J(\theta,u)$ in the product
space $M=\R\times E$, in which we introduce a metric $\norm\cdot_{(\theta,u)}$ by \ref[4.1].
We set for $F\in T_{(\theta,u)}^*M$
    $$  \norm F_{(\theta,u),*} = \sup\{ F(\kappa,v);\, 
        (\kappa,v)\in T_{(\theta,u)}M, \ \norm{(\kappa,v)}_{(\theta,u)} \leq 1\}.
    $$
The standard distance $\dist_M$ on $M$ is given by
    $$  \eqalign{
        &\dist_M((\theta_0,u_0),(\theta_1,u_1))  \cr
        &= \inf\left\{ \int_0^1\norm{\dot\sigma(t)}_{\sigma(t)}\, dt;\,
        \sigma\in C^1([0,1],M),\, \sigma(i)=(\theta_i,u_i)\ \hbox{for}\ i=0,1\right\}.  \cr}
    $$
Writing
    $$  D=(\partial_\theta,d_u),
    $$
we have
    $$  DJ(\theta,u)[(\kappa,v)]
        = \partial_\theta J(\theta,u)\kappa + d_uJ(\theta,u)v. 
    $$
By the definition \ref[4.2] of $J$, 
    $$  \eqalign{
        &J(\theta,u) = J(0,\Phi_\theta u) = I(\Phi_\theta u), \cr
        &\partial_\theta J(\theta,u)=\partial_\theta J(0,\Phi_\theta u)= P(\Phi_\theta u), \cr
        &d_u J(\theta,u)v = d_u I(\Phi_\theta u)\Phi_\theta v
            \quad \hbox{for}\ v\in T_u S. \cr}
    $$
Thus
    $$  \norm{DJ(\theta,u)}_{(\theta,u),*} 
        = \left( \abs{P(\Phi_\theta u)}^2 + \norm{dI(\Phi_\theta u)}_{T_{\Phi_\theta u}^*S}^2\right)^{1/2}
                                                            \eqno\label[4.6]
    $$
We note that for $(\theta+\alpha,u)$, $(\theta_0+\alpha,u_0)$, $(\theta_1+\alpha,u_1)\in M$
    $$  \eqalignno{
        &\norm{DJ(\theta+\alpha,u)}_{(\theta+\alpha,u),*}
        = \norm{DJ(\theta,\Phi_\alpha u)}_{(\theta,\Phi_\alpha u),*},   \cr
        &\dist_M((\theta_0+\alpha,u_0), (\theta_1+\alpha,u_1))
        =\dist_M((\theta_0,\Phi_\alpha u_0), (\theta_1,\Phi_\alpha u_1)). &\label[4.7]\cr} 
    $$
For $b\in\R$, we set
    $$  \wK_b=\{ (\theta,u)\in M;\ J(\theta,u)=b, \, DJ(\theta,u)=0\}.
    $$
We note that $\wK_b$ is invariant under $\Phi_\theta$ and 
    $$  \hbox{$(\theta,u)\in \wK_b$ if and only if $\Phi_\theta u\in \widehat K_b$ for $\theta\in\R$.}   
    $$
By \ref[4.7]
    $$  \eqalignno{
        \dist_M((\theta,u),\wK_b) &= \dist_M((0,\Phi_\theta u),\wK_b) &\label[4.8]\cr
        &\leq \dist_S(\Phi_\theta u, \widehat K_b)           &\label[4.9] \cr}
    $$
Here $\dist_S(\cdot,\cdot)$ is the standard distance on $S$, that is, for $u_0$, $u_1\in S$
    $$  \dist_S(u_0,u_1) = \inf\left\{\int_0^1\norm{\dot\gamma(t)}_E;\,
        \gamma(t)\in C^1([0,1],S),\, \gamma(i)=u_i \, \hbox{for}\, i=0,1\right\}.
    $$

\medskip

\proclaim Lemma \label[Lemma:4.7].
Assume that $I(u)$ satisfies $(PSP)_b$ on $S$.  Then
\item{(i)} Let $\{ (\theta_n,u_n)\}_{n=1}^\infty\subset M$ be a $(PS)$ sequence for $J$ at level $b$,
that is,
    $$  J(\theta_n,u_n)\to b \quad \hbox{and}\quad
        \norm{DJ(\theta_n,u_n)}_{(\theta_n,u_n),*}\to 0 \quad \hbox{as}\ n\to \infty.
    $$
Then $\{ \Phi_{\theta_n}u_n\}_{n=1}^\infty$ has a strongly convergent subsequence in $S$.
Moreover $\widehat K_b\not=\emptyset$ and 
    $$  \dist_M((\theta_n,u_n),\wK_b) \to 0 \quad \hbox{as}\ n\to \infty.   \eqno\label[4.10]
    $$
\item{(ii)} Suppose $\widehat K_b\not=\emptyset$, equivalently $\wK_b\not=\emptyset$.  Then for any $\rho>0$
there exists $\delta_\rho>0$ such that
    $$  \norm{DJ(\theta,u)}_{(\theta,u),*} \geq \delta_\rho
    $$
if $J(\theta,u)\in [b-\delta_\rho,b+\delta_\rho]$ and $(\theta,u)\not\in \wN_\rho(\wK_b)$.  Here
    $$  \wN_\rho(\wK_b)=\{ (\theta,u)\in M;\, \dist_M((\theta,u),\wK_b) <\rho\}.
    $$
\item{(iii)} 
If $\widehat K_b=\emptyset$, equivalently $\wK_b=\emptyset$, there exists $\delta_0>0$ such that
    $$  \norm{DJ(\theta,u)}_{(\theta,u),*} \geq \delta_0
    $$
for $(\theta,u)\in M$ with $J(\theta,u)\in [b-\delta_0,b+\delta_0]$.

\medskip

\claim Proof.
(i) Suppose that $\{ (\theta_n,u_n)\}_{n=1}^\infty$ is a $(PS)_b$ sequence for $J$ at level $b$.
By \ref[4.6], $\{ \Phi_{\theta_n} u_n\}_{n=1}^\infty$ satisfies $I(\Phi_{\theta_n}u_n)\to b$,
$P(\Phi_{\theta_n}u_n)\to b$, $\norm{dI(\Phi_{\theta_n}u_n)}_{T_{\Phi_{\theta_n}u_n}^*S}\to 0$.
Since $I(u)$ satisfies $(PSP)_b$ condition, (i) follows. \m
Moreover by \ref[4.9] we have \ref[4.10] and (ii), (iii) follow easily from (i).  \QED

\medskip

Following Palais [\cite[P:25]], we have

\medskip

\proclaim Corollary \label[Corollary:4.8].  
Set $\wK=\{ (\theta,u)\in M;\, DJ(\theta,u)=0\}$.  Then there exists a locally Lipschitz vector
field $W\in \calX(M\setminus\wK)$ such that for $(\theta,u)\in M\setminus\wK$
    $$  \eqalign{
        &\norm{W(\theta,u)}_{(\theta,u)} \leq 2\norm{DJ(\theta,u)}_{(\theta,u),*}, \cr
        &DJ(\theta,u)W(\theta,u) \geq \norm{DJ(\theta,u)}_{(\theta,u),*}^2. \cr}
    $$
Moreover under \ref[4.4]--\ref[4.5]
    $$  W_1(\theta,-u) = W_1(\theta,u), \quad W_2(\theta,-u) = -W_2(\theta,u),
    $$
where we write $W(\theta,u)=(W_1(\theta,u), W_2(\theta,u))\in T_{(\theta,u)}M$.

\medskip

We consider the following ODE in $M$
    $$  \left\{\eqalign{
        &{d\weta\over dt} = -\varphi(\weta)\psi(J(\weta)) {W(\weta)\over \norm{W(\weta)}_{\weta}}, \cr
        &\weta(0,\theta,u)=(\theta,u), \cr}
        \right.             \eqno\label[4.11]
    $$
where $\varphi:\, M\to [0,1]$, $\psi:\,\R\to [0,1]$ are locally Lipschitz cut-off functions such
that for small $\rho>0$
    $$  \eqalign{
    &\varphi(\theta,u)=\cases{  1   &for $(\theta,u)\not\in M\setminus\wN_{{2\over 3}\rho}(\wK_b)$, \cr
                                0   &for $(\theta,u)\in \wN_{{1\over 3}\rho}(\wK_b)$, \cr} \cr
    &\psi(t)=\cases{    1   &for $t\in [b-{\overline\epsilon\over 2},b+{\overline\epsilon\over 2}]$, \cr
                        0   &for $t\not\in [b-\overline\epsilon,b+\overline\epsilon]$. \cr} \cr}
    $$
We note that under \ref[4.4]--\ref[4.5], we may assume $\varphi(\theta,-u)=\varphi(\theta,u)$.

To show our Proposition \ref[Proposition:4.5], we need the following lemma, in which we use 
notation for $c\in \R$
    $$  [J\leq c]_M = \{ (\theta,u)\in M;\, J(\theta,u)\leq c\}.
    $$
By \ref[4.11], we may prove

\medskip

\proclaim Lemma \label[Lemma:4.9].
Suppose that $\overline\epsilon>0$ and $\rho>0$.  Then there exist $\epsilon\in (0,\overline\epsilon)$
and $\weta\in C([0,1]\times M,M)$ such that
\item{(i)} $\weta(0,\theta,u)=(\theta,u)$ for $(\theta,u)\in M$.
\item{(ii)} $\weta(t,\theta,u)=(\theta,u)$ for $t\in [0,1]$ 
if $(\theta,u)\in [J\leq b-\overline\epsilon]_M$.  
\item{(iii)} $t\mapsto J(\weta(t,\theta,u))$ is non-increasing for $(\theta,u)\in M$.
\item{(iv)} $\weta(1, [J\leq b+\epsilon]_M\setminus \wN_\rho(\wK_b)) \subset [J\leq b-\epsilon]_M$,
$\weta(1, [J\leq b+\epsilon]_M) \subset [J\leq b-\epsilon]_M\cup \wN_\rho(\wK_b))$. \m
When $\wK_b=\emptyset$, equivalently $\widehat K_b=\emptyset$, we regard $\wN_\rho(\wK_b)=\emptyset$.

{\sl\noindent
Moreover, if $S$ is symmetric with respect to $0$ and $I(u)$ is even in $u$,
\item{(v)} $\weta(t,\theta,u)=(\weta_1(t,\theta,u),\weta_2(t,\theta,u))$ satisfies
    $$  \weta_1(t,\theta,-u)=\weta_1(t,\theta,u), \quad 
        \weta_2(t,\theta,-u)=-\weta_2(t,\theta,u).   \QED
    $$
}

To deduce our Proposition \ref[Proposition:4.5], we need the following operators
    $$  \eqalign{
        &\iota:\, S\to M; \, u\mapsto (0,u), \cr
        &\pi:\, M\to S;\, (\theta,u)\mapsto \Phi_\theta u.\cr}
    $$
We have

\medskip

\proclaim Lemma \label[Lemma:4.10].
For any $\rho>0$ there exists an $R(\rho)>0$ such that
    $$  \eqalignno{
        &\pi(\wN_\rho(\wK_b)) \subset N_{R(\rho)}(\widehat K_b),                         &\label[4.12]\cr
        &\iota(S\setminus N_{R(\rho)}(\widehat K_b)) \subset M\setminus \wN_\rho(\wK_b), &\label[4.13]\cr
        &R(\rho)\to 0 \quad \hbox{as}\ \rho\to 0,                               &\label[4.14]\cr}
    $$
where
    $$  N_r(\widehat K_b)=\{ u\in S;\, \dist_S(u,\widehat K_b) <r \}.
    $$

\claim Proof.
First we show \ref[4.12].  Suppose that $(\theta,u)\in \wN_\rho(\wK_b)$.  By \ref[4.8], note that
$\dist_M((0,\Phi_\theta u),\wK_b)=\dist_M((\theta,u),\wK_b)<\rho$ and choose a 
$\sigma(t)\in C^1([0,1],M)$ such that $\sigma(0)=(0,\Phi_\theta u)$, $\sigma(1)\in \wK_b$,
$\int_0^1\norm{{d\sigma\over dt}(t)}_{\sigma(t)} \, dt<\rho$.  
Writing $\sigma(t)=(\sigma_1(t),\sigma_2(t))$, we have
    $$  \abs{\sigma_1(t)} \leq \int_0^1 \abs{{d\sigma\over dt}(t)}\, dt 
        \leq \int_0^1 \norm{{d\sigma\over dt}(t)}_{\sigma(t)}\, dt  <\rho.
    $$
We note that there exists $c_\rho>0$ such that for some $\delta\in (0,1]$
    $$  \eqalign{
		&c_\rho \norm v_E \leq \norm{\Phi_\theta v}_E \quad \hbox{for}\ \abs\theta\leq\rho \ 
			\hbox{and}\ v\in E, \cr
		&\delta\leq c_\rho \leq 1 \quad \hbox{for}\ \rho\in (0,1].\cr}
    $$
Thus
    $$  \eqalign{
        \dist_S(\Phi_\theta u,\sigma_2(1)) 
        &\leq \int_0^1\norm{{d\sigma_2\over dt}(t)}_E \, dt
        \leq c_\rho^{-1}\int_0^1 \norm{\Phi_{\sigma_1(t)} {d\sigma_2\over dt}(t)}_E \, dt \cr
        &\leq c_\rho^{-1}\int_0^1 \norm{{d\sigma_2\over dt}(t)}_{\sigma(t)} \, dt
        \leq c_\rho^{-1}\rho. \cr}
    $$
Therefore
    $$  \eqalign{
        \dist_S(\pi(\theta,u),\widehat K_b) &= \dist_S(\Phi_\theta u,\widehat K_b)
        \leq \dist_S(\Phi_\theta u,\Phi_{\sigma_1(1)}\sigma_2(1)) \cr
        &\leq \dist_S(\Phi_\theta u,\sigma_2(1)) 
            + \dist_S(\sigma_2(1), \Phi_{\sigma_1(1)}\sigma_2(1)) \cr
        &\leq c_\rho^{-1}\rho +\sup\{\dist_S(w,\Phi_\alpha w); 
            \abs{\alpha}\leq\rho,\, w\in \widehat K_b\}.\cr}
    $$
Since $\widehat K_b$ is compact by $(PSP)_b$, we have
    $$  R(\rho)=c_\rho^{-1}\rho +\sup\{\dist_S(w,\Phi_\alpha w); 
            \abs{\alpha}\leq\rho,\, w\in \widehat K_b\} \to 0 \quad \hbox{as}\ \rho\to 0
    $$
and \ref[4.12] and \ref[4.14] hold.

On the other hand, if $u\in S\setminus N_{R(\rho)}(\widehat K_b)$, by \ref[4.12] we have $\iota(u)=(0,u)
\in M\setminus\wN_\rho(\wK_b)$.  Thus \ref[4.13] holds.  \QED

\medskip

\claim Proof of Proposition \ref[Proposition:4.5].  
For a given neighborhood $O$ of $\widehat K_b$, first we choose $\rho>0$ so small that $N_{R(\rho)}(\widehat K_b)\subset O$.
By Lemma \ref[Lemma:4.10], we have $\pi(\wN_\rho(\wK_b))\subset N_{R(\rho)}(\widehat K_b)$.

For any $\overline\epsilon>0$, there exist $\epsilon\in (0,\overline\epsilon)$ and 
$\weta\in C([0,1]\times M,M)$ with the properties stated in Lemma \ref[Lemma:4.9].  Then we define 
    $$  \eta(t,u)=\pi(\weta(t,0,u)):\, [0,1]\times S\to S.
    $$
Properties (i)--(iii), (v) in Proposition \ref[Proposition:4.5] are easily checked.
As to (iv), we note that
    $$  \iota([I\leq b+\epsilon]_S\setminus O) 
        \subset \iota ([I\leq b+\epsilon]_S\setminus N_{R(\rho)}(\widehat K_b))
        \subset [J\leq b+\epsilon]_M\setminus \wN_\rho(\wK_b),
    $$
from which we have
    $$  \eqalign{
        \eta([I\leq b+\epsilon]_S\setminus O) 
        &= \pi(\weta(1,\iota([I\leq b+\epsilon]_S\setminus O))
        \subset \pi(\weta(1, [J\leq b+\epsilon]_M\setminus \wN_\rho(\wK_b))) \cr
        &\subset \pi([J\leq b-\epsilon]_M) \subset [I\leq b-\epsilon]_S.\cr}
    $$
Similarly, we have
    $$  \eqalign{
        \eta(1,[I\leq b+\epsilon]_S) 
        &\subset \pi(\weta(1,\iota([I\leq b+\epsilon]_S))) 
        \subset \pi(\weta(1,[J\leq b+\epsilon]_M)) \cr
        &\subset \pi([J\leq b-\epsilon]_M\cup \wN_{\rho}(\wK_b))
        \subset [I\leq b-\epsilon]_S\cup N_{R(\rho)}(\widehat K_b)  \cr
        &\subset [I\leq b-\epsilon]_S\cup O.\cr}
    $$
\QED

\bigskip


\bibliography

\medskip

\parindent=1.5\parindent

\bibitem[AKY:1] T. Akahori, H. Kikuchi, T. Yamada, 
Virial functional and dynamics for nonlinear Schr\"odinger equations of local interactions, 
NoDEA Nonlinear Differential Equations Appl. 25 (2018), no. 1, Art. 5, 27.

\bibitem[AN:2] T. Akahori, H. Nawa, 
Blowup and scattering problems for the nonlinear Schr\"odinger equations, 
Kyoto J. Math. 53 (2013), no. 3, 629--672.

\bibitem[AdAP:3] A. Azzollini, P. d'Avenia, A. Pomponio, 
Multiple critical points for a class of nonlinear functionals, Ann. Mat. Pura
Appl. (4) 190 (3) (2011) 507--523.

\bibitem[AP:4] A. Azzollini, A. Pomponio, 
On the Schr\"odinger equation in $\R^N$ under the effect of a general nonlinear term, 
Indiana Univ. Math. J. 58 (2009), no. 3, 1361--1378.

\bibitem[BL:5] A. Bahri, P.-L. Lions,
On the existence of a positive solution of semilinear elliptic
equations in unbounded domains,
Ann. Inst. H. Poincar\'e Anal. Non Lineair\'e 14 (1997), no. 3,
365--413 (1997).

\bibitem[BdV:6] T. Bartsch, S. de Valeriola, 
Normalized solutions of nonlinear Schr\"odinger equations,
Arch. Math. (Basel) 100 (2013), no. 1, 75--83.

\bibitem[BS1:7] T. Bartsch, N. Soave,
A natural constraint approach to normalized solutions of nonlinear Schr\"odinger equations and systems,
J. Funct. Anal. 272 (2017), no. 12, 4998--5037 and 
Correction, J. Funct. Anal. 275 (2018), no. 2, 516--521.

\bibitem[BS2:8] T. Bartsch, N. Soave, 
Multiple normalized solutions for a competing system of Schr\"odinger equations, 
Calc. Var. Partial Differential Equations 58 (2019), no. 1, 58:22.

\bibitem[BeL1:9] H. Berestycki, P.-L. Lions,
Nonlinear scalar field equations. I. Existence of a ground state,
Arch. Rational Mech. Anal. 82 (1983), no. 4, 313--345.

\bibitem[BeL2:10] H. Berestycki, P.-L. Lions,
Nonlinear scalar field equations. I. Existence of infinitely many solutions,
Arch. Rational Mech. Anal. 82 (1983), no. 4, 347--375.

\bibitem[BeGK:11] H. Berestycki, T. Gallou\"et, O. Kavian,
\'Equations de champs scalaires euclidiens non lin\'eaires dans le plan,
C. R. Acad. Sci. Paris Ser. I Math. 297 (1983), no. 5, 307--310.

\bibitem[BT:12] J. Byeon, K. Tanaka, 
Nonlinear Elliptic Equations in Strip-Like Domains,
Advanced Nonlinear Studies 12 (2012), 749--765.

\bibitem[CT:13] C.-N. Chen, K. Tanaka, 
A variational approach for standing waves of FitzHugh-Nagumo type systems,
J. Diff. Eq. 257 (2014) 109--144.

\bibitem[FISJ:14] G. M. Figueiredo, N. Ikoma, J. R. Santos Junior, 
Existence and concentration result for the Kirchhoff type equations with general nonlinearities, 
Arch. Ration. Mech. Anal. 213 (2014), no. 3, 931--979.

\bibitem[G:15] N. Ghoussoub,
Duality and perturbation methods in critical point theory, 
Cambridge Tracts in Mathematics 107, Cambridge University Press, Cambridge (1993).

\bibitem[HT:16] J. Hirata, K. Tanaka,
Nonlinear scalar field equations with $L^2$ constraint: Mountain pass and symmetric mountain pass approaches,
Advances Nonlinear Studies (to appear).

\bibitem[HIT:17] J. Hirata, N. Ikoma, K. Tanaka,
Nonlinear scalar field equations in {$\R^N$}: mountain pass and symmetric mountain pass approaches,
Topol. Methods Nonlinear Anal. 35 (2010), no. 2, 253--276.

\bibitem[I1:18] N. Ikoma, 
Existence of solutions of scalar field equations with fractional operator, 
J. Fixed Point Theory Appl. 19 (2017), no. 1, 649--690 and 
Erratum, J. Fixed Point Theory Appl. 19 (2017), no. 2, 1649--1652.

\bibitem[I2:19] N. Ikoma,
Multiplicity of radial and nonradial solutions to equations with fractional operators,
preprint (2019).

\bibitem[J:20] L. Jeanjean,
Existence of solutions with prescribed norm for semilinear elliptic equations,
Nonlinear Anal. 28 (1997), no. 10, 1633--1659.

\bibitem[JL:21] L. Jeanjean, S.-S. Lu, 
Nonradial normalized solutions for nonlinear scalar field equations, Preprint (2018). 

\bibitem[LM:22] R. Lehrer, L. A. Maia, 
Positive solutions of asymptotically linear equations via Pohozaev manifold, 
J. Funct. Anal. 266 (2014), no. 1, 213--246.

\bibitem[MMP:23] R. Mandel, E. Montefusco, B. Pellacci, 
Oscillating solutions for nonlinear Helmholtz equations, 
Z. Angew. Math. Phys. 68 (2017), no. 6, Art. 121, 19pp.

\bibitem[MVS:24] V. Moroz, J. Van Schaftingen, 
Existence of ground states for a class of nonlinear Choquard equations,
Trans. Amer. Math. Soc. 367, no. 9, 6557--6579.

\bibitem[P:25] R. S. Palais, 
Lusternik-Schnirelman theory on Banach manifolds
Topology 5 (1966), 115--132.

\bibitem[R:26] P. H. Rabinowitz,
Minimax methods in critical point theory with applications to differential equations,
CBMS Regional Conference Series in Mathematics, 65 (1986).

\bibitem[Sh:27] J. Shatah, 
Unstable ground state of nonlinear Klein-Gordon equations, 
Trans. Amer. Math. Soc. 290 (1985), no. 2, 701--710.

\bibitem[Si:28] M. Shibata,
Stable standing waves of nonlinear Schr\"odinger equations with a general nonlinear term,
Manuscripta Math. 143 (2014), no. 1-2, 221--237.


\bye


\input date-file-info.mac


\vfil
\break

{

\baselineskip 28 pt


\def\cmd#1#2{\line{\indent\hbox to 2.5 cm{\tt \string#1\hfil}\hskip 1 cm $\displaystyle #1$\hfil}}

\let\def=\cmd

\centerline{\bf --- macros ---}

\bigskip

\input macro.mac

}


\bye